\documentclass[11pt,reqno]{amsart}

\usepackage[margin=1.15in]{geometry}
\usepackage{amssymb,mathtools,array,multirow}
\usepackage{microtype}
\usepackage{graphicx,calc}
\graphicspath{{figures/}}
\usepackage[caption=false,font=footnotesize,margin=0pt,justification=centering]{subfig}
\usepackage{tikz}
\usepackage{placeins,flafter,float}
\usepackage[round,authoryear]{natbib}
\usepackage[hidelinks]{hyperref}
\usepackage{bookmark}
\AddToHook{cmd/section/before}{\FloatBarrier}
\AddToHook{cmd/subsection/before}{\FloatBarrier}
\newcommand{\panelpair}[3][.485]{%
  \noindent\makebox[.485\linewidth]{\includegraphics[width=#1\linewidth,clip]{#2}}\hfill%
  \makebox[.485\linewidth]{\includegraphics[width=#1\linewidth,clip]{#3}}%
  \par\vspace{2pt}%
}

\newcommand{\R}{{\mathbb{R}}}
\newcommand{\Z}{{\mathbb{Z}}}
\newcommand{\F}{{\mathbb{F}}}
\newcommand{\cF}{\mathcal{F}}
\newcommand{\cX}{\mathcal{X}}
\newcommand{\cM}{\mathcal{M}}
\newcommand{\cZ}{\mathcal{Z}}
\newcommand{\olF}{\overline{F}}
\newcommand{\old}{\overline{d}}
\newcommand{\olf}{\overline{f}}
\newcommand{\olr}{\overline{r}}
\newcommand{\olcF}{\overline{\cF}}
\newcommand{\olgamma}{\overline{\gamma}}
\newcommand{\mvmap}{{\rightrightarrows}}
\newcommand{\setof}[1]{\left\{ {#1}\right\}}
\newcommand{\setdef}[2]{\left\{ #1 \, \middle| \, #2 \right\}}
\newcommand{\Int}{\mathop{\mathrm{int}}\nolimits}
\newcommand{\cl}{\mathop{\mathrm{cl}}\nolimits}
\newcommand{\Inv}{\mathop{\mathrm{Inv}}\nolimits}
\theoremstyle{plain}
\newtheorem{thm}{Theorem}[section]
\newtheorem{lem}[thm]{Lemma}
\newtheorem{prop}[thm]{Proposition}
\newtheorem{cor}[thm]{Corollary}
\theoremstyle{definition}
\newtheorem{defn}[thm]{Definition}
\theoremstyle{remark}
\newtheorem{rem}[thm]{Remark}

\title{An Equivariant Approach to Topological Dynamics}
\author[B. Rivas]{Bernardo Rivas}
\address{University of Toledo, Toledo, OH}
\email{bernardo.dopradorivas@utoledo.edu}
\author[K. Iwasaki]{Kaito Iwasaki}
\address{University of Michigan, Ann Arbor, MI}
\email{kaitoi@umich.edu}
\author[R. Vasudevan]{Ram Vasudevan}
\address{University of Michigan, Ann Arbor, MI}
\email{ramv@umich.edu}
\makeatletter
\hypersetup{pdftitle={\@title},pdfauthor={\authors},pdfkeywords={\@keywords}}
\makeatother

\begin{document}
\raggedbottom

\begin{abstract}%
Computational topological methods for analyzing global dynamics suffer from the curse of dimensionality and, in efficient implementations, rely on restrictive topological conditions. These limitations make the methods difficult to apply to high-dimensional manifolds, including the state spaces of robotic systems. We extend a framework for latent dynamics to equivariant systems by taking the orbit space as the latent space. We learn the encoder and latent map from data using equivariant networks. Equivariant outer approximations then reduce the combinatorial computation to the latent grid. A bound on the latent prediction error modulo the group certifies that latent attracting blocks pullback to attracting blocks for the original dynamics. Under additional hypotheses, we reconstruct Conley indices in the original space from lifted index pairs. We illustrate the framework with learned models of an $\mathrm{SO}(2)$-equivariant flow on the sphere and an $\mathrm{O}(3)$-equivariant flow, and with a quadrotor model learned from full rigid-body trajectories under reduced attitude control.%
\end{abstract}

\maketitle

\section{Introduction}
\label{sec:introduction}
Homological dynamics uses finite cell complexes and directed graphs to characterize global dynamics through Conley-theoretical frameworks~\citep{kalies:mischaikow:vandervorst:26}. For a continuous map, a cubical grid supplies a cell complex, and an outer approximation determines a directed graph by enclosing the images of the grid cells~\citep{arai:kalies:kokubu:mischaikow:oka:pilarczyk,bush:gameiro:harker:kokubu:mischaikow:obayashi:pilarczyk}. The cost of these computations, however, grows rapidly with dimension, and implementations based on rectangular grids and cubical homology~\citep{kaczynski:mischaikow:mrozek:04} restrict the state spaces that can be treated directly. These limitations motivate learning reduced models, especially for systems on high-dimensional manifolds such as those arising in robotics, that retrieve accurate and certified information.

Learned latent models have been used to represent nonlinear dynamics on low-dimensional manifolds~\citep{linot:graham:20} and to analyze robotic controllers with Morse graphs~\citep{morals}. In a Conley-theoretical setting, it has been shown that neural approximations support rigorous conclusions~\citep{gameiro:gelb:mischaikow}, while the latent-dynamics framework certifies the pullback of attracting blocks through a learned encoder under a bound on the dynamical residual~\citep{latent:dynamics:26}.

For equivariant systems, symmetry reduction~\citep{budanur:cvitanovic:davidchack:siminos:15} and symmetry-aware autoencoders~\citep{kneer:sayadi:sipp:schmid:rigas:21} provide reduced descriptions based on the group action. We extend latent dynamics to this setting by performing the combinatorial computation on the orbit space. We identify outer approximations of the quotient dynamics $\olF=F/G$ on $X/G$ with those of an equivariant map $F$ on the corresponding grids of $G$-invariant cells, so that the computation on $X/G$ preserves the Morse graph.

For learned equivariant models, we certify the pullback of latent attracting blocks using a bound on the prediction residual modulo $G$. Since the quotient does not in general determine Conley indices in the original space, we reconstruct homology and index maps for rotations of $S^2$ and for two specified lifts of the quotient dynamics under the orthogonal action on $\R\times\R^3$ (Theorems~\ref{thm:sphere-reconstruction} and~\ref{thm:orthogonal-reconstruction}). We illustrate the framework on three dynamical systems using models learned with equivariant multilayer perceptrons (EMLPs)~\citep{finzi:welling:wilson:21}.

\section{Combinatorial and latent dynamics}
\label{sec:combinatorial}

Throughout this manuscript, $X$ is a compact metric space, $F\colon X\to X$ is continuous and $\F$ is a field. Unless stated otherwise, homology has coefficients in $\F$. For $N\subset X$, $\Inv(N,F)$ is the set of points through which a full trajectory in $N$ exists. A compact set $N$ is an isolating neighborhood if $\Inv(N,F)\subset\Int N$, and $\Inv(N,F)$ is then an isolated invariant set~\citep{conley:cbms}. A compact set $N$ is an attracting block if $F(N)\subset\Int N$.

For $x\in X$ and a full trajectory $\gamma=(\gamma_n)_{n\in\Z}$, let $\omega(x)=\bigcap_{n\geq0}\cl\setdef{F^k(x)}{k\geq n}$ and $\alpha(\gamma)=\bigcap_{n\geq0}\cl\setdef{\gamma_{-k}}{k\geq n}$. A Morse representation of $F$ is a finite poset $(\mathsf M,\leq)$ of pairwise disjoint nonempty compact sets $M\subset X$ with $F(M)=M$ such that every $\omega(x)$ is contained in some $M\in\mathsf M$, and every full trajectory $\gamma$ with $\gamma_0\notin\bigcup_{M\in\mathsf M}M$ has $\omega(\gamma_0)\subset M$ and $\alpha(\gamma)\subset M'$ for some $M<M'$ in $\mathsf M$~\citep{kalies:mischaikow:vandervorst:22}.

A pair of compact sets $L\subset N$ is an index pair for $S=\Inv(\cl(N\setminus L),F)$ if $S\subset\Int(N\setminus L)$, $F(L)\cap N\subset L$ and $F(N\setminus L)\subset N$. On the pointed quotient $N/L$, the map $F_{(N,L)}$ induced by $F$ sends $[x]$ to $[F(x)]$ for $x\in N\setminus L$ and fixes the basepoint. When $L=\emptyset$, we take $N/L=N\sqcup\{*\}$, so $\widetilde H_*(N/L;\F)\cong H_*(N;\F)$. The index map is the endomorphism induced by $F_{(N,L)}$ on $\widetilde H_*(N/L;\F)$, and the Conley index of $S$ is its shift equivalence class~\citep{szymczak:95,franks:richeson:00}. A pair $(N,L)$ is good if $L$ is nonempty, closed and a deformation retract of a neighborhood in $N$. For a good pair the quotient map gives $\widetilde H_*(N/L;\F)\cong H_*(N,L;\F)$~\citep[Proposition~2.22]{hatcher:02}.

A grid of $X$ is a finite set $\cX$ of cells $\xi$ with nonempty compact regular-closed geometric realizations $|\xi|\subset X$ that have pairwise disjoint interiors and union $X$. For $\mathcal A\subset\cX$, let $|\mathcal A|=\bigcup_{\xi\in\mathcal A}|\xi|$. A multivalued map $\cF\colon\cX\mvmap\cX$ acts on subsets by $\cF(\mathcal A)=\bigcup_{\xi\in\mathcal A}\cF(\xi)$, and it is an outer approximation of $F$ if
\begin{equation}
  \label{eq:outer_approximation}
  F(|\xi|)\subset\Int|\cF(\xi)|\quad\forall\xi\in\cX,
\end{equation}
in which case $F$ is a selector of $\cF$.

Regarding $\cF$ as a directed graph on $\cX$, with an edge from $\xi$ to each $\eta\in\cF(\xi)$, a Morse set of $\cF$ is a strongly connected component $\cM\subset\cX$ that contains an edge. We use the name for both $\cM$ and $|\cM|$, distinguishing them by context or notation. The Morse sets are ordered by $\cM\leq\cM'$ if a path leads from $\cM'$ to $\cM$, and the Morse graph of $\cF$ is the Hasse diagram of this order. A set $\mathcal A\subset\cX$ is forward invariant if $\cF(\mathcal A)\subset\mathcal A$ and an attractor of $\cF$ if $\cF(\mathcal A)=\mathcal A$. The forward invariant sets form a finite distributive lattice under union and intersection, and the attractors form a finite distributive lattice under inclusion~\citep{kalies:mischaikow:vandervorst:14,kalies:mischaikow:vandervorst:16,kalies:mischaikow:vandervorst:22}.

\begin{figure}[H]
  \centering
  \begin{minipage}[b]{0.35\textwidth}\centering
    \subfloat[Outer approximation]{\includegraphics[scale=0.88]{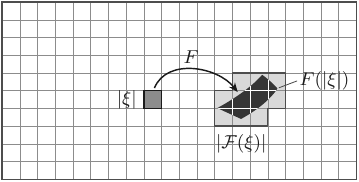}}
  \end{minipage}\hfill
  \begin{minipage}[b]{0.35\textwidth}\centering
    \subfloat[Morse sets]{\includegraphics[scale=0.88]{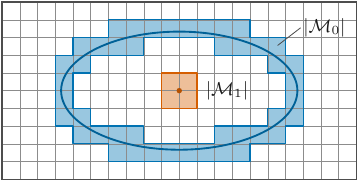}}
  \end{minipage}\hfill
  \begin{minipage}[b]{0.28\textwidth}\centering
    \subfloat[Conley--Morse graph]{\parbox[b][\heightof{\includegraphics[scale=0.88]{outer-approximation.pdf}}][c]{\linewidth}{\centering\includegraphics{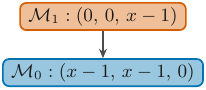}}}
  \end{minipage}
  \caption{A map with an unstable fixed point and an attracting invariant curve. In (c), each node lists, by degree, the characteristic polynomial of the index map, with $0$ for a zero group.}
  \label{fig:combinatorial}
\end{figure}

Let $\cF$ be an outer approximation of $F$. If $\mathcal A$ is forward invariant, then $F(|\mathcal A|)\subset\Int|\mathcal A|$, so $|\mathcal A|$ is an attracting block. For each Morse set $\cM$, the region $|\cM|$ is an isolating neighborhood. The nonempty sets $\Inv(|\cM|,F)$ form a Morse representation of $F$ with the order of the Morse graph~\citep{kalies:mischaikow:vandervorst:22}. For each Morse set $\cM$, the pair $\left(|\cM\cup\cF(\cM)|,|\cF(\cM)\setminus\cM|\right)$ is an index pair for $\Inv(|\cM|,F)$. On a cubical grid, \texttt{CMGDB}~\citep{CMGDB} computes the index map from a cubical enclosure of $F$ when the enclosure satisfies the acyclicity conditions of the relative homology algorithm on this index pair~\citep{arai:kalies:kokubu:mischaikow:oka:pilarczyk,bush:gameiro:harker:kokubu:mischaikow:obayashi:pilarczyk}. The Morse graph with the Conley index of $\Inv(|\cM|,F)$ at the node of each Morse set $\cM$ is the Conley--Morse graph (Figure~\ref{fig:combinatorial}).

When $E\colon X\to Z$ is a continuous encoder and $f\colon Z\to Z$ a continuous latent map on the compact metric space $(Z,d_Z)$, the tolerance of an attracting block $N\subset Z$ for $f$ is defined by $\tau_f(N)=\sup\setdef{\nu>0}{B_\nu(f(N))\subset N}\in(0,\infty]$, where $B_\nu(A)=\setdef{z \in Z}{d_Z(z,A)<\nu}$. The following theorem allows the attracting blocks to be pulled back when the residual $r(x)=d_Z\left(f(E(x)),E(F(x))\right)$ is sufficiently small.
\begin{thm}[\citealt{latent:dynamics:26}]
  \label{thm:main_simple}
  If $N$ is an attracting block for $f$ and $\sup_{x\in E^{-1}(N)}r(x)<\tau_f(N)$, then $E^{-1}(N)$ is an attracting block for $F$.
\end{thm}

\section{Equivariant outer approximations}
\label{sec:equivariant-outer}

Let $G$ be a compact Hausdorff group acting continuously on $X$, and assume that $F$ is equivariant, that is, $gF(x)=F(gx)$ for every $g \in G$ and $x \in X$. Let $\pi\colon X\to X/G$ be the orbit map and $\olF=F/G$ the induced map, and recall that the orbit space $X/G$ is compact metrizable.

\begin{defn}
\label{def:equivariant-grid}
A grid $\cX$ of $X$ is equivariant if for all $g\in G$ and $\xi\in\cX$ there is $\xi'\in\cX$ with $|\xi'|=g|\xi|$. Denote $g\xi\coloneqq \xi'$. A multivalued map $\cF\colon\cX\mvmap\cX$ on an equivariant grid $\cX$ is equivariant if $g\cF(\xi)\subset\cF(g\xi)$ for all $g\in G$ and $\xi\in\cX$.
\end{defn}
Since realizations have disjoint interiors, $\xi'$ is unique, so $\xi\mapsto g\xi$ is an action of $G$ on $\cX$, and the inclusion for $g^{-1}$ yields $\cF(g\xi)=g\cF(\xi)$. Thus $\cF$ induces a multivalued map on $\cX/G$ (Proposition~\ref{prop:quotient-outer}), and when every $|\xi|$ is $G$-invariant, the quotient grid has the same cells $\cX$, with realizations $\pi(|\xi|)$.

\begin{prop}
\label{prop:pullback-grid}
Let $\cX$ be a finite set and let $|\xi|\subset X$ be $G$-invariant for every $\xi\in\cX$.
\begin{enumerate}
\item[(i)] $\cX$ is a grid of $X$ with realizations $|\xi|$ if and only if it is a grid of $X/G$ with realizations $\pi(|\xi|)$. Moreover, every grid of $X/G$ is of this form.
\item[(ii)] Let $\cX$ be a grid of $X$ and $F\colon X\to X$ be equivariant. Every multivalued map $\cF\colon\cX\mvmap\cX$ is equivariant, and it is an outer approximation of $F$ if and only if it is one of $\olF$ on the quotient grid.
\item[(iii)] For $\cF$ as in (ii), its Morse sets, forward invariant sets and attractors are the same subsets of $\cX$ for both grids. The orders on the Morse sets, and hence the Morse graphs, are equal. For each $\mathcal A\subset\cX$, the corresponding regions are $|\mathcal A|\subset X$ and $\pi(|\mathcal A|)\subset X/G$.
\end{enumerate}
\end{prop}
\begin{proof}
The orbit map $\pi$ is continuous, open, closed and surjective. Continuity and openness imply, for every $B\subset X/G$,
\[
  \pi^{-1}(\Int B)=\Int\pi^{-1}(B)\quad\text{and}\quad\pi^{-1}(\cl B)=\cl\pi^{-1}(B).
\]
For (i), each $G$-invariant realization satisfies $|\xi|=\pi^{-1}(\pi(|\xi|))$. The identities above and injectivity of inverse image imply that regular closedness and disjointness of interiors hold for the realizations in $X$ if and only if they hold for their images in $X/G$. Closed realizations are compact in both spaces, and surjectivity ensures equivalence of nonemptiness and covering. The preimages of the realizations of any grid on $X/G$ are $G$-invariant, so every grid of $X/G$ is of this form.

For (ii), the action $g\xi=\xi$ satisfies $|g\xi|=g|\xi|$ and $g\cF(\xi)=\cF(g\xi)$ for every multivalued map $\cF$. Since $|\cF(\xi)|$ is $G$-invariant, the interior identity implies that $F(|\xi|)\subset\Int|\cF(\xi)|$ if and only if $\pi(F(|\xi|))\subset\Int\pi(|\cF(\xi)|)$. The latter is the outer approximation condition for $\olF$, because $\pi\circ F=\olF\circ\pi$.

For (iii), both grids use the same directed graph of $\cF$ on $\cX$, so the Morse sets and their order, forward invariant sets, attractors and Morse graphs coincide. Since $\pi$ commutes with unions, the region of $\mathcal A\subset\cX$ on the quotient grid is $\bigcup_{\xi\in\mathcal A}\pi(|\xi|)=\pi(|\mathcal A|)$, which proves (iii).
\end{proof}

\begin{prop}
\label{prop:invariant-sets}
Let $F\colon X\to X$ be equivariant.
\begin{enumerate}
\item[(i)] If $K\subset X$ is $G$-invariant, then so is $\Inv(K,F)$, and $\Inv(K,F)/G=\Inv(K/G,\olF)$.
\item[(ii)] Let $N$ and $L$ be $G$-invariant subsets of $X$. If $(N/G,L/G)$ is an index pair for $\olF$, then $(N,L)$ is an index pair for $F$, and $S=\Inv(\cl(N\setminus L),F)$ is $G$-invariant with $S/G=\Inv(\cl(N/G\setminus L/G),\olF)$.
\end{enumerate}
\end{prop}
\begin{proof}
For (i), equivariance implies $F(Gx)=GF(x)$, hence
\begin{equation}
  \label{eq:preimage_commutes}
  F(\pi^{-1}(w))=\pi^{-1}(\olF(w))\quad\forall w\in X/G.
\end{equation}
Let $\olgamma$ be a full trajectory of $\olF$ in $K/G$ and let $x\in\pi^{-1}(\olgamma_0)$. Define $\gamma_n=F^n(x)$ for $n\geq0$ and recursively choose $\gamma_{-n-1}\in\pi^{-1}(\olgamma_{-n-1})$ with $F(\gamma_{-n-1})=\gamma_{-n}$ for $n\geq0$, whose existence follows from \eqref{eq:preimage_commutes}. Then $\pi(\gamma_n)=\olgamma_n$ for every $n\in\Z$, so $\gamma$ is a full trajectory in $K=\pi^{-1}(K/G)$ through $x$. Together with projection of full trajectories, this proves that $\Inv(K,F)=\pi^{-1}(\Inv(K/G,\olF))$. This set is $G$-invariant, and taking its image under $\pi$ proves (i).

For (ii), the sets $N=\pi^{-1}(N/G)$ and $L=\pi^{-1}(L/G)$ are compact in $X$, since $N/G$ and $L/G$ are compact, and $L\subset N$. The identity $\pi\circ F=\olF\circ\pi$ and the index-pair conditions on $X/G$ imply $F(L)\cap N\subset L$ and $F(N\setminus L)\subset N$. Since preimages under the continuous open map $\pi$ commute with interiors and closures, (i) and the isolation condition on $X/G$ imply
\begin{align*}
  S=\pi^{-1}(\Inv(\cl((N/G)\setminus(L/G)),\olF)) \subset\pi^{-1}(\Int((N/G)\setminus(L/G)))=\Int(N\setminus L).
\end{align*}
Thus $(N,L)$ is an index pair for $F$, and $S$ is $G$-invariant with $S/G=\Inv(\cl((N/G)\setminus(L/G)),\olF)$.
\end{proof}

By Proposition~\ref{prop:pullback-grid}, an outer approximation $\cF$ of $\olF$ on the quotient of a grid of $X$ with $G$-invariant realizations is also one of $F$, so the conclusions of Section~\ref{sec:combinatorial} hold in both spaces. The invariant sets $\Inv(|\cM|,F)$ and the index pairs of the Morse sets in $X$ are $G$-invariant and project to those for $\olF$, by Proposition~\ref{prop:invariant-sets}(i) and since $\pi$ commutes with unions. In particular, a region $|\cM|$ that contains $S\subset X$ contains every $gS$ with $g\in G$.

\begin{rem}
\label{rem:sign}
Note that the Conley indices of a $G$-invariant isolated invariant set and its quotient can differ, and equivariant maps with the same quotient dynamics can have different Conley indices. For example, let $G=\mathrm{SO}(3)$ act on $X=S^2$ and let $F_\pm(x)=\pm x$ be equivariant maps. Note that they induce the same map on $X/G=\setof{*}$, the pair $(X,\emptyset)$ is an index pair for both maps, and their index maps are the identity on $H_0(X)=\F$ and multiplication by $\pm1$ on $H_2(X)=\F$. On $X/G$, the index map of both is simply the identity on $\F$ in degree $0$ and does not preserve degree $2$. Conversely, the index on $X/G$ does not determine the index on $X$ since the indices differ when $\operatorname{char}\F\neq2$.
\end{rem}

\section{Residual bounds and attracting blocks}
\label{sec:learning}

Let a compact group $G$ act continuously on $X$ and on the compact metric space $(Z,d_Z)$, with the action on $Z$ by isometries. Let $F\colon X\to X$ and $f\colon Z\to Z$ be continuous and equivariant, and let $E\colon X\to Z$ be a continuous equivariant encoder. If $E$ is surjective and $E\circ F=f\circ E$, then $E$ is a semiconjugacy from $F$ to $f$. When the action on $Z$ is trivial, equivariance of $E$ means $E(gx)=E(x)$, so $E$ factors through the orbit map $\pi\colon X\to X/G$ as $E=\widetilde E\circ\pi$. If $E(x)=E(y)$ is equivalent to $Gx=Gy$ for all $x,y\in X$, then $\widetilde E$ is a homeomorphism from $X/G$ onto $E(X)$. Moreover, if $E(X)=Z$, then $Z \cong X/G$.

For the action on $Z$, let $\pi_Z\colon Z\to Z/G$ be the orbit map, $\olf=f/G$, and $A/G=\pi_Z(A)$ for a $G$-invariant $A\subset Z$. Recall that the metric $\old(\pi_Z(a),\pi_Z(b))=\min_{g\in G}d_Z(a,gb)$ induces the quotient topology on $Z/G$.

\begin{defn}
\label{def:residual-modulo-G}
The residual modulo $G$ is the map $\olr\colon X\to[0,\infty)$ defined for $x\in X$ by
\[
  \olr(x)=\min_{g\in G}d_Z\left(f(E(x)),gE(F(x))\right)=\old\left(\olf(\pi_Z(E(x))),\pi_Z(E(F(x)))\right).
\]
\end{defn}
Thus $\olr$ measures the failure of the semiconjugacy identity after projection to $Z/G$ and satisfies $\olr(x)\leq r(x)$ for the residual $r$ of Section~\ref{sec:combinatorial}. The equivariance and isometry properties imply that $r(kx)=r(x)$ and $\olr(kx)=\olr(x)$ for $k\in G$, using the substitution $g\mapsto k^{-1}gk$. Thus a residual bound on a $G$-invariant region can be checked on a subset meeting every orbit in that region.

Let $\cZ$ be a grid of $Z$ with $G$-invariant realizations, and let $\cF_Z\colon\cZ\mvmap\cZ$ be an outer approximation of $\olf$ on the quotient grid, with realizations $\pi_Z(|\zeta|)$ for $\zeta\in\cZ$ (Proposition~\ref{prop:pullback-grid}).

\begin{prop}
\label{prop:equivariant_pullback}
Let $\mathcal A\subset\cZ$ with $\cF_Z(\mathcal A)\subset\mathcal A$ and let $N=|\mathcal A|$.
\begin{enumerate}
\item[(i)] $N/G$ is an attracting block for $\olf$, $N$ is a $G$-invariant attracting block for $f$, and $\tau_f(N)=\tau_{\olf}(N/G)$.
\item[(ii)] Let $S\subset E^{-1}(N)$ meet every $G$-orbit in $E^{-1}(N)$. If $\sup_{x\in S}\olr(x)<\tau_f(N)=\tau_{\olf}(N/G)$, then $E^{-1}(N)$ is a $G$-invariant attracting block for $F$.
\item[(iii)] The set $\omega(E^{-1}(N),F)=\bigcap_{n\geq0}\cl\bigcup_{k\geq n}F^k(E^{-1}(N))$ is $G$-invariant.
\end{enumerate}
\end{prop}
\begin{proof}
For (i), forward invariance of $\mathcal A$ and the outer approximation property imply $\olf(N/G)\subset\Int(N/G)$. Proposition~\ref{prop:pullback-grid}(ii) implies $f(N)\subset\Int N$. For every $G$-invariant $A\subset Z$ and $\nu>0$, the quotient metric satisfies
\[
  \pi_Z^{-1}\left(B_\nu(A/G)\right)=B_\nu(A).
\]
Applying this identity to $f(N)$ shows that $B_\nu(f(N))\subset N$ if and only if $B_\nu(\olf(N/G))\subset N/G$, so the tolerances are equal.

For (ii), constancy on orbits implies that $\sup_{x\in E^{-1}(N)}\olr(x)=\sup_{x\in S}\olr(x)$. Thus, Theorem~\ref{thm:main_simple}, applied to $\pi_Z\circ E$ and $\olf$, yields that $(\pi_Z\circ E)^{-1}(N/G)=E^{-1}(N)$ is an attracting block for $F$. Finally, equivariance of $E$ and invariance of $N$ imply that $E^{-1}(N)$ is $G$-invariant.

For (iii), equivariance of $F$ makes every forward image of $E^{-1}(N)$ invariant, thus their unions, closures and intersections are invariant, and the same holds for $\omega(E^{-1}(N),F)$.
\end{proof}

The forward invariant sets of cells satisfying the residual bound are closed under intersections, and finite unions and intersections of their realized regions pull back to $G$-invariant attracting blocks for $F$ (Corollary~\ref{prop:block-lattice}).

\section{Examples}
\label{sec:examples}

We compute Conley--Morse graphs with an adapted version of \texttt{CMGDB}~\citep{CMGDB}, using uniform dyadic grids.\footnote{Code: \url{https://github.com/bernardorivas/2027_L4DC_EquivariantMorseGraphs}.} For a homeomorphism $\Theta\colon X/G\to B$ onto a box $B\subset\R^k$, a product grid with $m=2^\ell$ intervals per axis has $m^k=2^{\ell k}$ cells $K_\xi$. By Proposition~\ref{prop:pullback-grid}, the realizations $|\xi|=(\Theta\circ\pi)^{-1}(K_\xi)$ form a $G$-invariant grid of $X$, and outer approximations of $F$ on this grid coincide with those of $\Theta\circ\olF\circ\Theta^{-1}$. The attracting blocks, isolating neighborhoods and Morse representations obtained from an outer approximation are valid for every selector of that map.

We parameterize the learned maps by EMLPs~\citep{finzi:welling:wilson:21} whose representations enforce the required invariance or equivariance, label each reference node by the invariant set of $F$ in its region, and, since the computed graphs agree, assign each learned node the label of the corresponding reference node.

Let $\cM$ be a minimal Morse set of the outer approximation of a learned map, or a union of Morse sets merged into one node. The set $\mathcal A$ of cells from which no path leaves the box or reaches a Morse set outside $\cM$ is forward invariant, and we call $N=|\mathcal A|$ the region of attraction of $\cM$. For the inequality of Proposition~\ref{prop:equivariant_pullback}(ii), we estimate the tolerance of $N$ by the distance from the realization of the image of $\mathcal A$ under the outer approximation to the complement of $N$, and the supremum of the residual modulo $G$ by its maximum over a sample of $E^{-1}(N)$, both for the map whose Morse graph is computed.

The reference indices in $X$ are proved on explicit $G$-invariant index pairs, whose index maps are the identity by the homotopy along each flow (Appendix~\ref{supp:example-pairs}), and the correspondence between learned nodes and these invariant sets is observed.

\subsection{An \texorpdfstring{$\mathrm{SO}(2)$}{SO(2)}-equivariant flow on \texorpdfstring{$S^2$}{S2}}
\label{sec:ex:so2}

For $p = (x,y,z) \in X = S^2$, consider the vector field $V(p)=0.6 e_3 \times p+z(0.4^2-z^2)(e_3-z p)$, where $e_3 = (0,0,1)$. The group $G=\mathrm{SO}(2)$ acts by rotating $(x,y)$ and leaving $z$ fixed, and the vector field and its flow are equivariant under this action. Let $F$ be the time-ten map of this flow.

Each orbit is at a specific latitude, so the coordinate $z$ identifies the quotient $X/G$ with $[-1,1]$. Under this identification, $\olF$ is the time-ten map of $\dot z=z(0.4^2-z^2)(1-z^2)$. In the quotient, the poles $\pm e_3$ and the equator $C=\setof{z=0}$ correspond to repelling fixed points, while the latitudes $A_\pm=\setof{z=\pm0.4}$ correspond to attracting fixed points.

Following this setup, we jointly fit an invariant encoder $E\colon S^2\to[-1,1]$, a learned time-one map $f\colon[-1,1]\to[-1,1]$ with $f(\pm1)=\pm1$ by construction, and a decoder $D\colon[-1,1]\to S^2$ that reconstructs a representative of each rotation orbit. We iterate $f$ ten times to obtain $f^{10}$, then compare its Conley--Morse graph with that of $\olF$ in the height coordinate.

\begin{figure}[H]
\centering
\begin{minipage}[b]{.36\linewidth}\centering
\subfloat[Predicted height after ten steps]{\makebox[\linewidth]{\includegraphics[height=1.35in]{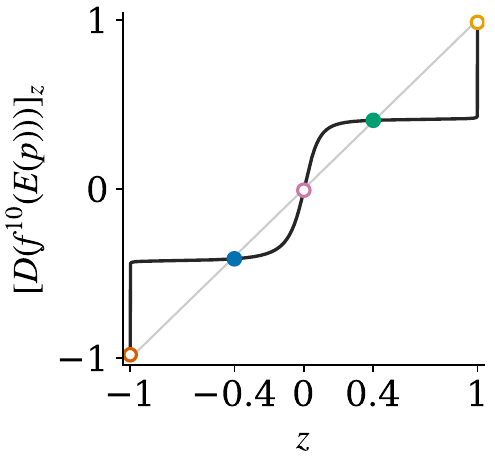}}}
\end{minipage}\hfill
\begin{minipage}[b]{.62\linewidth}\centering
\subfloat[{Conley--Morse graph in $[-1,1]$}]{\parbox[b][1.35in][c]{\linewidth}{\centering\includegraphics[width=.9\linewidth]{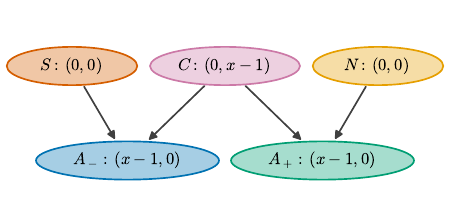}}}
\end{minipage}
\caption{Ten-step predictions (a) and the Conley--Morse graph of $f^{10}$ (b).}
\label{fig:sphere-computation}
\end{figure}

On a grid of $2^7$ cells, the learned computation has six Morse sets. Two adjacent Morse sets representing the learned counterpart of $C$ are combined, and CMGDB recomputes the index on the merged index pair. The resulting five-node graph has the same directed connections and indices as the reference graph (Figure~\ref{fig:appendix-sphere-comparison} in Appendix~\ref{supp:sphere-walkthrough}). For the learned blocks $N$ associated with $A_-$ and $A_+$, the maximum sampled residual modulo $G$ of $f^{10}$ on each block is smaller than its estimated tolerance, as required by Proposition~\ref{prop:equivariant_pullback}(ii), supporting their pullback as attracting blocks $E^{-1}(N)$ for $F$, whose attractors $\omega(E^{-1}(N),F)$ are $G$-invariant by Proposition~\ref{prop:equivariant_pullback}(iii) (Appendix~\ref{supp:sphere-computation}).

\begin{figure}[H]
\centering
\begin{minipage}[b]{.36\linewidth}\centering
\subfloat[$E^{-1}(|\cM|)$ in $S^2$]{\makebox[\linewidth]{\includegraphics[height=1.62in]{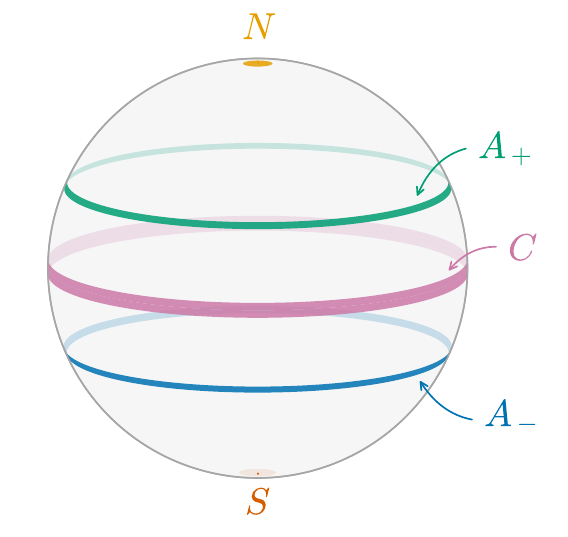}}}
\end{minipage}\hfill
\begin{minipage}[b]{.62\linewidth}\centering
\subfloat[Reconstructed Conley--Morse graph in $S^2$]{\parbox[b][1.62in][c]{\linewidth}{\centering\includegraphics[width=.9\linewidth]{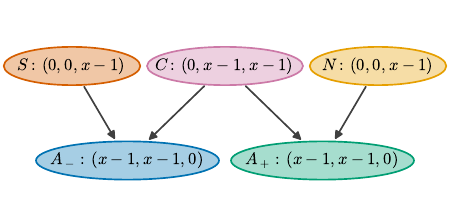}}}
\end{minipage}
\caption{Preimages of the Morse sets of $f^{10}$ and the reconstructed indices in $S^2$.}
\label{fig:sphere-reconstruction}
\end{figure}

Applying Theorem~\ref{thm:sphere-reconstruction} to the computed index pairs of $f^{10}$ yields the indices of its equivariant lifts to $S^2$. The indices of the reconstructed Conley--Morse graph in Figure~\ref{fig:sphere-reconstruction}(b) agree with the analytical indices of $F$ (Appendix~\ref{supp:sphere-pairs} and Figure~\ref{fig:appendix-sphere-comparison}).

\subsection{An \texorpdfstring{$\mathrm{O}(3)$}{O(3)}-equivariant flow on \texorpdfstring{$\R^4$}{R4}}
\label{sec:ex:o3}

For a scalar $a\in\R$ and a vector $v\in\R^3$, consider the system $\dot a=a(1-a^2)$, $\dot v=(1-2a^2-\|v\|^2)v$. The group $G=\mathrm{O}(3)$ acts on $(a,v)$ by orthogonal transformations of $v$, leaving $a$ fixed, and the vector field and its flow are equivariant under this action. We restrict the flow to the forward-invariant set $X=[-3/2,3/2]\times B^3_{3/2}$, where $B^3_{3/2}=\setof{v\in\R^3:\|v\|\leq3/2}$. Let $F$ denote the time-$\log3$ map of this flow. Each group orbit consists of states with the same scalar $a$ and norm $\|v\|$, so these coordinates identify $X/G$ with $[-3/2,3/2]\times[0,3/2]$.

\begin{figure}[H]
\centering
\begin{minipage}[b]{.48\linewidth}\centering
\subfloat[Exact quotient flow]{\makebox[\linewidth]{\includegraphics[height=.52\linewidth]{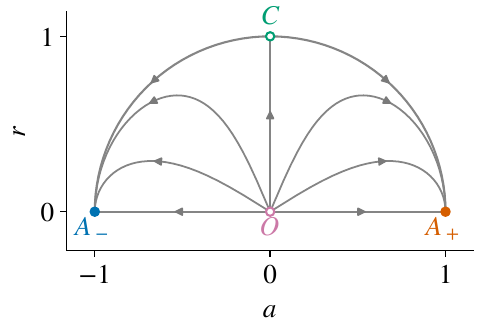}}}
\end{minipage}\hfill
\begin{minipage}[b]{.48\linewidth}\centering
\subfloat[Morse sets and regions of attraction]{\makebox[\linewidth]{\includegraphics[width=.7533\linewidth,height=.52\linewidth]{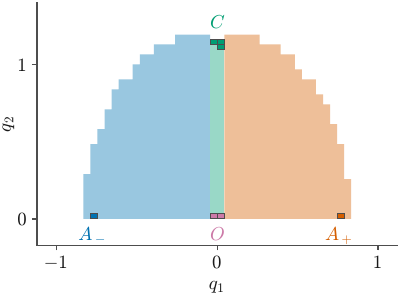}}}
\end{minipage}

\begin{minipage}[b]{.46\linewidth}\centering
\subfloat[Learned Conley--Morse graph in $H$]{\makebox[\linewidth]{\includegraphics[scale=.52,trim=0 0 0 19,clip]{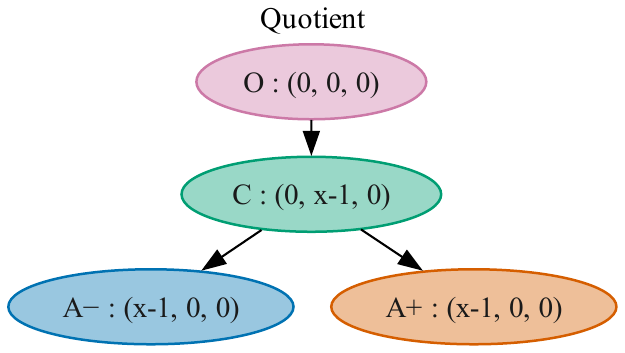}}}
\end{minipage}\hfill
\begin{minipage}[b]{.50\linewidth}\centering
\subfloat[Reconstructed Conley--Morse graph in $X$]{\makebox[\linewidth]{\includegraphics[scale=.52,trim=0 0 0 22,clip]{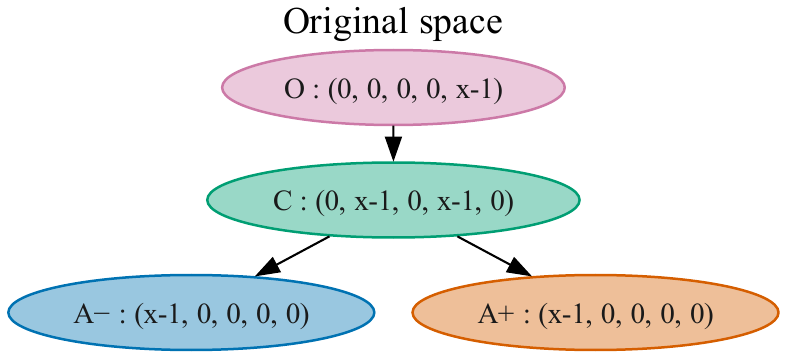}}}
\end{minipage}
\caption{The $\mathrm O(3)$-equivariant flow and its learned model.}
\label{fig:o3}
\end{figure}

In the quotient flow, the origin $O=(0,0)$ is a repelling fixed point, the sphere of equilibria $C=\setof{a=0,\ \|v\|=1}$ projects to the saddle fixed point $(0,1)$, and $A_\pm=(\pm1,0)$ are attracting fixed points. There are connecting trajectories from $O$ to the saddle and from the saddle to each attractor, and the invariant semicircle $a^2+\|v\|^2=1$ contains the connections from the saddle to the two attractors, as shown in Figure~\ref{fig:o3}(a).

Let $H=\setof{(q_1,q_2)\in\R^2:q_2\geq0}$. We jointly fit an invariant encoder $E\colon X\to H$, a learned time-$\log3$ map $f\colon H\to H$, and a decoder $D\colon H\to\R^4$ that reconstructs a representative of each $\mathrm O(3)$-orbit. We compare the Conley--Morse graph of $f$ with that of the quotient map $\olF$. We compute the learned graph on a rectangular grid in $H$, using bounding boxes of sampled images under $f$ to define the multivalued map supplied to CMGDB.

On a $2^6\times2^6$ grid, the learned computation has six Morse sets. We combine the two Morse sets associated with $O$ and the two associated with $C$, including the cells on directed paths between the Morse sets in each pair. CMGDB then recomputes the indices on the merged index pairs, yielding a four-node graph (Appendix~\ref{supp:o3-grouping}). The resulting graph has the same directed connections and indices as the reference quotient graph (Figure~\ref{fig:appendix-o3-comparison} in Appendix~\ref{supp:o3-walkthrough}).

To reconstruct indices in $X$, we lift the learned map $f$ by preserving the direction of the vector coordinate and apply Theorem~\ref{thm:orthogonal-reconstruction} to local quotient pairs, which are separate from the merged CMGDB cell pairs. The indices of the resulting graph in Figure~\ref{fig:o3}(d) agree with the analytical indices of $F$ (Appendix~\ref{supp:o3-pairs} and Figure~\ref{fig:appendix-o3-comparison}).

\subsection{Reduced attitude control on \texorpdfstring{$S^2$}{S2}}
\label{sec:ex:attitude}

We next consider attitude control of a quadrotor. In the rigid-body model~\citep{lee:leok:mcclamroch:10,kousik:holmes:vasudevan:19}, the state is $(p,v,R,\Omega)\in\R^3\times\R^3\times\mathrm{SO}(3)\times\R^3$, where $p$ is position, $v$ is velocity, $R$ is orientation, and $\Omega$ is body angular velocity. The thrust direction $q=Re_3\in S^2$ determines the direction of the propulsive force and therefore affects translational motion. The control objective is to stabilize this direction at the upright pole $e_3$. Write $z=q\cdot e_3$ for its height.

For a smooth gain $a\colon[-1,1]\to\R$ and constant angular speed $\omega$, prescribe the thrust-axis vector field $V_a(q)=a(z)(e_3-zq)+\omega e_3\times q$. Let $w=R\Omega$ be the spatial angular velocity, $w_d(q)=q\times V_a(q)$ its desired value, and $\eta=w-w_d(q)$ the angular-velocity error. To make the angular velocity track $w_d(q)$, we apply constant thrust $\mu g$ and control moment
\[
  M=\Omega\times J\Omega+JR^\top\left(Dw_d(q)[w\times q]-2\eta\right),
\]
where $\mu$ is the mass, $g$ is gravitational acceleration, and $J$ is the body inertia matrix. The resulting dynamics satisfy $\dot q=V_a(q)+\eta\times q$ and $\dot\eta=-2\eta$. Thus the angular-velocity error decays exponentially, and on the invariant set $\setof{\eta=0}$, the thrust direction follows $\dot q=V_a(q)$.

The group $G=\mathrm{SO}(2)$ acts by rotations $Q$ about the vertical axis $e_3$, with $(p,v,R,\Omega)\mapsto(Qp,Qv,QR,\Omega)$ and $q\mapsto Qq$. The feedback and the thrust-axis dynamics are equivariant under this action. Let $\omega=\pi/\log2$, and let $F\colon X\to X$ be the time-one map of $V_a$ on $X=S^2$. As in Section~\ref{sec:ex:so2}, the height $z$ identifies $X/G$ with $[-1,1]$, where the quotient dynamics satisfy $\dot z=a(z)(1-z^2)$.

For the nominal gain $a=1$, the downward thrust direction $-e_3$ is a repelling fixed point, while the upright direction $e_3$ attracts every thrust direction except $-e_3$. For the coning gain $a(z)=1/4-z$, both poles repel, while every nonpolar thrust direction approaches the periodic orbit $C=\setof{z=1/4}$, rotating about the vertical axis. The orbit has period $2\log2$ and is locally exponentially orbitally stable in the reduced state $(q,\eta)$.

We test whether the learned Conley--Morse graphs distinguish upright stabilization from convergence to an unwanted coning orbit using a one-dimensional computation. Using plant trajectories initialized on $\setof{\eta=0}$, we jointly fit an equivariant encoder $E\colon\R^3\times\R^3\times\mathrm{SO}(3)\times\R^3\to Z$ into a latent sphere $Z$, an equivariant time-$0.1$ map $f\colon Z\to Z$, and a decoder $D\colon Z\to S^2$ that reconstructs the thrust direction. Training also penalizes the discrepancy between $E(p,v,R,\Omega)$ and the known direction $q=Re_3$. Using the height projection $\pi_Z\colon Z\to[-1,1]$, we identify $Z/G$ with $[-1,1]$ and learn an increasing homeomorphism $\psi\colon[-1,1]\to[-1,1]$. The resulting interval encoder is $\psi\circ\pi_Z\circ E$.

\begin{figure}[H]
\centering
\begin{minipage}[b]{.36\linewidth}\centering
\subfloat[{Learned Conley--Morse graph}]{\parbox[b][1.1in][c]{\linewidth}{\centering\includegraphics[scale=.55]{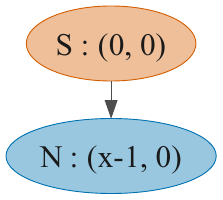}}}
\end{minipage}\hfill
\begin{minipage}[b]{.30\linewidth}\centering
\subfloat[Morse sets]{\makebox[\linewidth]{\includegraphics[height=1.1in]{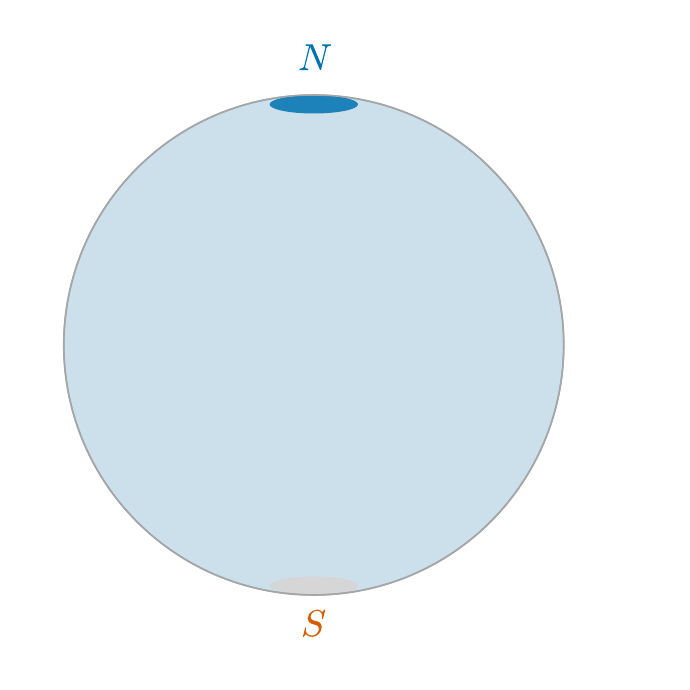}}}
\end{minipage}\hfill
\begin{minipage}[b]{.33\linewidth}\centering
\subfloat[Conley--Morse graph in $S^2$]{\parbox[b][1.1in][c]{\linewidth}{\centering\includegraphics[scale=.55]{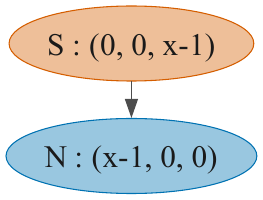}}}
\end{minipage}

\begin{minipage}[b]{.36\linewidth}\centering
\subfloat[{Learned Conley--Morse graph}]{\parbox[b][1.1in][c]{\linewidth}{\centering\includegraphics[scale=.55]{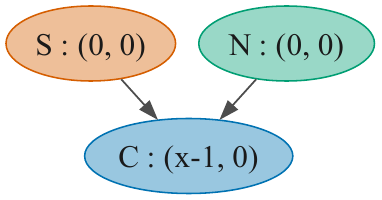}}}
\end{minipage}\hfill
\begin{minipage}[b]{.30\linewidth}\centering
\subfloat[Morse sets]{\makebox[\linewidth]{\includegraphics[height=1.1in]{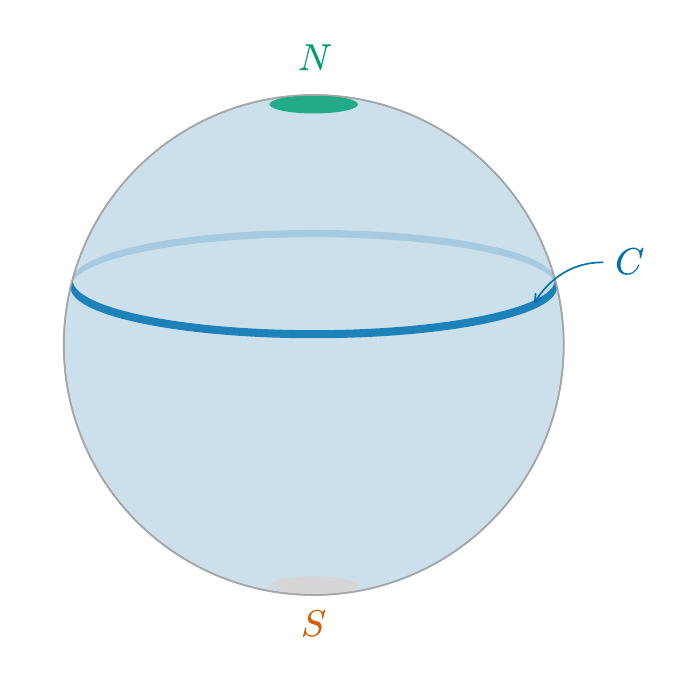}}}
\end{minipage}\hfill
\begin{minipage}[b]{.33\linewidth}\centering
\subfloat[Conley--Morse graph in $S^2$]{\parbox[b][1.1in][c]{\linewidth}{\centering\includegraphics[scale=.55]{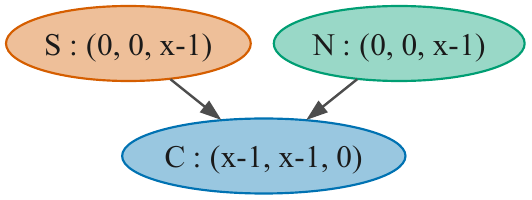}}}
\end{minipage}
\caption{Nominal gain (a)--(c) and coning gain (d)--(f).}
\label{fig:attitude-computation}
\end{figure}

We compute the Conley--Morse graph of $(\psi\circ\olf\circ\psi^{-1})^{10}$ on $2^7$ interval cells. For the nominal gain, the learned and reference quotient graphs have the same directed connection $-e_3\to e_3$ and interval Conley indices (Figure~\ref{fig:attitude-computation}(a)).

For the coning gain, the reference quotient graph in the height coordinate $z$ has directed connections $-e_3\to C\leftarrow e_3$. The learned computation yields an additional one-cell Morse set with zero Conley index adjacent to the Morse set labeled $C$. After merging the additional Morse set and recomputing its interval index, the learned graph has the reference order $-e_3\to C\leftarrow e_3$ (Figure~\ref{fig:attitude-computation}(d); Appendix~\ref{supp:attitude-grouping}). In this graph, the upright direction is repelling and the sink corresponds to $C$, identifying failure of upright stabilization. For both gains, the maximum sampled residual modulo $G$ of $f^{10}$ on $\setof{\eta=0}$ is smaller than the estimated tolerance of the attracting block, which for the coning gain is the block of the merged node (Appendix~\ref{supp:attitude-computation}).

Panels (c) and (f) of Figure~\ref{fig:attitude-computation} annotate the learned Morse graphs with the analytical Conley indices of the corresponding thrust-axis invariant sets (Corollary~\ref{cor:attitude-axis-index}). The index of $C$ is nonzero in degree one, which distinguishes the unwanted attracting circle from the desired attracting point, although both have the same quotient index.

\section{Discussion}
\label{sec:discussion}

The framework presented applies computational topology on orbit spaces and addresses two bottlenecks simultaneously: computational complexity and topological constraints. Factoring through the quotient allows the analysis of noncontractible spaces, and reduces the CMGDB time from $47$\,s to $3.1$\,ms in Section~\ref{sec:ex:so2}, from $1{,}011$\,s to $1.8$\,ms for Section~\ref{sec:ex:o3}. In Section~\ref{sec:ex:attitude}, the computation uses only $2^7$ cells, while a grid of the 12-dimensional state with same resolution would have $2^{84}$.

However, the benefits come with limitations. We highlight that retaining the group action is central to computing the Conley index under $G$~\citep{floer:87}, and that the quotient does not determine the Conley index even when the maps are known (Remark~\ref{rem:sign}), so we use the rotation action on $S^2$ in Theorem~\ref{thm:sphere-reconstruction} and the specified radial lift in Theorem~\ref{thm:orthogonal-reconstruction}. That said, the Morse graphs are computed directly, and in Section~\ref{sec:ex:attitude} they distinguish upright stabilization from coning: the nominal graph has an upright sink, while both poles are repelling nodes of the coning graph. For both systems, the maximum sampled residual is below the estimated tolerance of the attracting block (Appendix~\ref{supp:attitude-computation}).

A next step is to replace sampled data and residuals by validated enclosures and uniform bounds, so that the attracting-block criterion of Proposition~\ref{prop:equivariant_pullback}(ii) certifies conclusions for $F$. In the current setting, the attracting blocks $E^{-1}(N)$ for $F$ are observed rather than proved, and the agreement of the learned Conley--Morse graphs with the reference graphs is observed, not certified. Finally, the method assumes that $G$, its action and a homeomorphism $\Theta$ of $X/G$ onto a box are known and that $F$ is equivariant, so, unfortunately, unknown or approximate symmetries are outside its scope.

\section*{Generative AI Use}

Generative AI tools (Anthropic's Claude Code and OpenAI's Codex) were used for experimental design, mathematical derivations and proofs, coding and experimental implementation, interpretation of results, and writing and editing. The authors reviewed and verified all AI-assisted content and take responsibility for the final content of the paper.

\section*{Acknowledgments}

This work was supported by the Air Force Office of Scientific Research under award number FA9550-23-1-0400.

\bibliographystyle{plainnat}
\bibliography{references-arxiv}

\clearpage
\appendix
\section{Supporting proofs}
\label{app:supporting-proofs}

\subsection{Equivariant grids}
\label{proof:quotient-outer}

We use the setting of Section~\ref{sec:equivariant-outer}.

\begin{prop}
\label{prop:quotient-outer}
Let $\cX$ be an equivariant grid of $X$.
\begin{enumerate}
\item[(i)] $\cX/G$ is a grid of $X/G$ with realization $\pi(|\xi|)$ for each cell $G\xi$.
\item[(ii)] Let $F\colon X\to X$ be equivariant and let $\cF$ be an equivariant outer approximation of $F$ on $\cX$. Then $\olcF(G\xi)=\setdef{G\eta}{\eta\in\cF(\xi)}$ is well defined and is an outer approximation of $\olF$ on $\cX/G$.
\end{enumerate}
\end{prop}

\begin{proof}
Since $\pi(|g\xi|)=\pi(|\xi|)$, the sets $\pi(|\xi|)$ are indexed by $\cX/G$. The map $\pi$ is open and closed, so these sets are nonempty compact regular closed and cover $X/G$. If $G\xi\neq G\eta$, then
\[
  \pi^{-1}(\pi(|\xi|)\cap\pi(|\eta|))
  =\bigcup_{\substack{\eta'\in G\eta}}\left(\bigcup_{\substack{\xi'\in G\xi}}(|\xi'|\cap|\eta'|)\right).
\]
Note that each intersection on the right is closed with empty interior since $\xi'$ and $\eta'$ are distinct. Thus, the finite union has empty interior and the projected realizations have disjoint interiors, which proves (i). To prove (ii), note that equivariance implies $\cF(g\xi)=g\cF(\xi)$, so $\olcF$ is well defined. Since $|\olcF(G\xi)|=\pi(|\cF(\xi)|)$, the outer approximation property of $\cF$ and openness of $\pi$ imply
\begin{equation*}
  \olF(\pi(|\xi|))=\pi(F(|\xi|))\subset\pi(\Int|\cF(\xi)|)\subset\Int\pi(|\cF(\xi)|)=\Int|\olcF(G\xi)|.
\end{equation*}
\end{proof}

The minimal outer approximation of $F$ on a grid $\cX$ of $X$ satisfies~\eqref{eq:outer_approximation} and $\cF_{\min}(\xi)\subset\cF(\xi)$ for every outer approximation $\cF$ of $F$ and every $\xi\in\cX$. Explicitly, it is given by
\begin{equation}
  \label{eq:minimal_outer_approximation}
  \cF_{\min}(\xi)=\setdef{\eta\in\cX}{|\eta|\cap F(|\xi|)\neq\emptyset}\quad\forall\xi\in\cX.
\end{equation}
If the realizations of $\cX$ are $G$-invariant and $F$ is equivariant, then $F$ and $\olF$ have the same outer approximations on $\cX$ by Proposition~\ref{prop:pullback-grid}(ii), so their minimal outer approximations coincide.

\subsection{Lattices of attracting blocks}
\label{proof:block-lattice}

We use the notation and hypotheses of Section~\ref{sec:learning}.

\begin{cor}
\label{prop:block-lattice}
Let
\[
  \mathcal S=\setdef{\mathcal A\subset\cZ}{\cF_Z(\mathcal A)\subset\mathcal A,\ \sup_{x\in E^{-1}(|\mathcal A|)}\olr(x)<\tau_{\olf}(|\mathcal A|/G)}.
\]
Then $\mathcal S$ is closed under intersections, and $E^{-1}(N)$ is a $G$-invariant attracting block for $F$ for every $N$ in the lattice of subsets of $Z$ generated by $\setdef{|\mathcal A|}{\mathcal A\in\mathcal S}$.
\end{cor}

\begin{proof}
Let $\tau(\mathcal C)=\tau_{\olf}(|\mathcal C|/G)$ for forward invariant $\mathcal C\subset\cZ$. For $\mathcal A,\mathcal B\in\mathcal S$, the set $\mathcal A\cap\mathcal B$ is forward invariant. For $0<\nu<\min\{\tau(\mathcal A),\tau(\mathcal B)\}$, the definitions of the tolerances and the openness of the neighborhood imply
\[
  B_\nu\left(\olf(|\mathcal A\cap\mathcal B|/G)\right)\subset\Int(|\mathcal A|/G)\cap\Int(|\mathcal B|/G)=\Int(|\mathcal A\cap\mathcal B|/G),
\]
where the equality holds on the quotient grid. Therefore
\[
  \sup_{x\in E^{-1}(|\mathcal A\cap\mathcal B|)}\olr(x)<\min\{\tau(\mathcal A),\tau(\mathcal B)\}\leq\tau(\mathcal A\cap\mathcal B),
\]
so $\mathcal A\cap\mathcal B\in\mathcal S$ and Proposition~\ref{prop:equivariant_pullback}(ii) certifies every $E^{-1}(|\mathcal A|)$. Thus, closure of $G$-invariant attracting blocks under finite unions and intersections, and commutativity of these operations with $E^{-1}$, prove the lattice claim.
\end{proof}

\subsection{Homotopies of index maps}

For compact $L\subset N$ and a continuous map $T$ defined on $N$, the two exit inclusions are $T(L)\cap N\subset L$ and $T(N\setminus L)\subset N$. When they hold, write $T_{(N,L)}$ for the pointed map of Section~\ref{sec:combinatorial}.

\begin{lem}
\label{lem:pointed-homotopy}
Let $Y$ be a Hausdorff space and $L\subset N\subset Y$ compact sets, let $A$ be a compact Hausdorff space and $A_0\subset A$ a closed set, and let $h\colon A\times[0,1]\to Y$ be continuous with
\[
  h_s(A_0)\cap N\subset L\quad\text{and}\quad h_s(A\setminus A_0)\subset N\quad\forall s\in[0,1].
\]
Then the maps $\hat h_s\colon A/A_0\to N/L$ with $\hat h_s([a])=[h_s(a)]$ for $a\in A\setminus A_0$, which map the basepoint to the basepoint, form a homotopy of pointed maps.
\end{lem}

\begin{proof}
Let $q\colon N\to N/L$ be the canonical map and $*$ the basepoint of $N/L$. Define $\Psi(a,s)=q(h_s(a))$ if $h_s(a)\in N$ and $\Psi(a,s)=*$ otherwise, and let
\[
  A'=\setdef{(a,s)\in A\times[0,1]}{h_s(a)\in N},\quad C=\cl\setdef{(a,s)\in A\times[0,1]}{h_s(a)\notin N}.
\]
The set $A'$ is closed because $N$ is closed in $Y$, and $\Psi=q\circ h$ on $A'$. By the second inclusion $C\subset A_0\times[0,1]$, and for $(a,s)\in A_0\times[0,1]$ the first inclusion shows that $h_s(a)\notin N$ or $h_s(a)\in L$, so $\Psi=*$ on $C$. Since $A'$ and $C$ are closed and cover $A\times[0,1]$, $\Psi$ is continuous. It is constant on $A_0\times[0,1]$, and $\Psi(a,s)=\hat h_s([a])$ for $a\in A\setminus A_0$. If $A_0=\emptyset$, then $A/A_0$ is $A$ with a disjoint basepoint, and $\hat h$ is continuous. Otherwise $A\times[0,1]\to(A/A_0)\times[0,1]$ is a continuous surjection from a compact space onto a Hausdorff space, hence a quotient map, and $\hat h$ is continuous.
\end{proof}

With $A=N$, $A_0=L$ and $h$ constant in $s$, Lemma~\ref{lem:pointed-homotopy} shows that $T_{(N,L)}$ is continuous whenever $T$ is continuous and the two inclusions hold. If $h_0$ is the inclusion of $N$, it shows that $(h_1)_{(N,L)}$ is homotopic to the identity, so it induces the identity on $\widetilde H_*(N/L)$.

\begin{lem}
\label{lem:preimage-pair}
Let $G$ act by homeomorphisms on a Hausdorff space $X$ with Hausdorff orbit space $X/G$, and let $\pi\colon X\to X/G$ be the orbit map. Let $\olF\colon X/G\to X/G$ be continuous and $L\subset N\subset X$ compact $G$-invariant sets with
\[
  \olF(L/G)\cap(N/G)\subset L/G\quad\text{and}\quad\olF((N\setminus L)/G)\subset N/G.
\]
Let $F_s\colon N\to X$, $s\in[0,1]$, be a homotopy with $\pi\circ F_s=\olF\circ\pi$ for all $s$. Then
\[
  F_s(L)\cap N\subset L\quad\text{and}\quad F_s(N\setminus L)\subset N\quad\forall s\in[0,1].
\]
If $(N/G,L/G)$ is an index pair for $\olF$, then $(N,L)$ is an index pair for every $F_s$. The pointed maps $(F_0)_{(N,L)}$ and $(F_1)_{(N,L)}$ are homotopic.
\end{lem}

\begin{proof}
Since $N$ and $L$ are $G$-invariant and $\pi\circ F_s=\olF\circ\pi$, for every $s\in[0,1]$ we have
\begin{align*}
  F_s(L)\cap N&\subset\pi^{-1}\left(\olF(L/G)\cap(N/G)\right)\subset L,\\
  F_s(N\setminus L)&\subset\pi^{-1}\left(\olF((N\setminus L)/G)\right)\subset N.
\end{align*}
Let $(N/G,L/G)$ be an index pair and $\gamma$ a full trajectory of $F_s$ in $\cl(N\setminus L)$. The continuous open map $\pi$ pulls closures and interiors back to the closures and interiors of the preimages. Thus $\pi\circ\gamma$ is a full trajectory in $\cl((N/G)\setminus(L/G))$. Isolation downstairs gives $\pi(\gamma(0))\in\Int((N/G)\setminus(L/G))$, hence $\gamma(0)\in\Int(N\setminus L)$. The last statement is Lemma~\ref{lem:pointed-homotopy} with $A=N$, $A_0=L$ and $h_s=F_s$.
\end{proof}

\section{Reconstruction of Conley indices}
\label{app:reconstruction}

One summand in the reconstruction uses the quotient index pair $(B,A)$. The other uses $(B,A\cup D)$, where $D$ consists of the points in $B$ whose fibers are singletons. These fibers occur at the pole heights in Sections~\ref{sec:ex:so2} and~\ref{sec:ex:attitude} and on the axis $r=0$ in Section~\ref{sec:ex:o3}.

For a continuous map $h$ and compact sets $C\subset B$ satisfying $h(C)\cap B\subset C$ and $h(B\setminus C)\subset B$, let $\psi_C\colon B/C\to B/C$ be the pointed map of Section~\ref{sec:combinatorial}. For the finite CW pairs below, $(\psi_C)_{*,j}$ denotes the induced endomorphism of $H_j(B,C)\cong\widetilde H_j(B/C)$, with $B/\emptyset=B\sqcup\{*\}$. We take $H_j=0$ for $j<0$.

\subsection{The \texorpdfstring{$\mathrm{SO}(2)$}{SO(2)} action on \texorpdfstring{$S^2$}{S2}}
\label{app:sphere-instance}

Let $\mathrm{SO}(2)$ act on $S^2$ by rotations about $e_3$. The height $\pi(w,z)=z$ identifies the quotient with $I=[-1,1]$. With $\rho(z)=\sqrt{1-z^2}$, the map $\Phi(z,u)=(\rho(z)u,z)$ parametrizes the latitude at $z$ by $u\in S^1$. At each endpoint of $I$ the whole circle represents one pole.

\begin{thm}
\label{thm:sphere-reconstruction}
Let $F\colon S^2\to S^2$ be continuous and $\mathrm{SO}(2)$-equivariant, and let $h\colon I\to I$ be its height map. Let $(B,A)$ be an index pair for $h$, with $B$ nonempty and $A$ and $B$ subcomplexes of a finite subdivision of $I$. Define $D=B\cap\{-1,1\}$, $K=A\cup D$, $N=\pi^{-1}(B)$ and $L=\pi^{-1}(A)$. Then $(N,L)$ is an index pair for $F$. In every degree $k\geq0$, the index map and the direct sum on the right are conjugate.
\begin{align*}
  J_k\colon \widetilde H_k(N/L)&\xrightarrow{\cong} H_k(B,A)\oplus H_{k-1}(B,K),\\
  J_k\,(F_{(N,L)})_{*,k}\,J_k^{-1}&=(\psi_A)_{*,k}\oplus(\psi_K)_{*,k-1}.
\end{align*}
\end{thm}

\begin{proof}
Equivariance sends each pole to a pole, so $h(D)\subset\{-1,1\}$. The exit inclusions for $(B,A)$ imply $h(K)\cap B\subset K$ and $h(B\setminus K)\subset B$, which define $\psi_K$. Proposition~\ref{prop:invariant-sets}(ii) lifts the index pair to $(N,L)$.

Let $o=(1,0)\in S^1$. The meridian $\sigma(z)=\Phi(z,o)$ embeds $B/A$ as a pointed subspace $M$ of $W=N/L$, with inverse the restriction of the height projection $\bar\pi\colon W\to B/A$. Give $S^1$ one vertex $o$ and one open $1$-cell. In $W$, the circle over each point of $D$ collapses to its meridian point, and all fibers over $A$ collapse to the basepoint. Besides the basepoint, the cells are $e\times\{o\}$ for open cells $e$ of $B$ outside $A$ and $e\times(S^1\setminus\{o\})$ for $e$ outside $K$. Their product characteristic maps descend to $W$, making $M$ a CW subcomplex.

Let $q\colon W\to W/M$ collapse the meridian. In $W/M$ all of $K\times S^1$ and $B\times\{o\}$ go to the basepoint, and all remaining points of $B\times S^1$ are distinct. The resulting continuous bijection of compact Hausdorff spaces identifies
\[
 W/M\cong (B/K)\wedge S^1=\Sigma(B/K).
\]
Since the section composed with $\bar\pi$ retracts $W$ onto $M$, the reduced exact sequence of $(W,M)$ splits. Let $\operatorname{susp}\colon H_{k-1}(B,K)\to\widetilde H_k(W/M)$ be the suspension isomorphism~\citep{hatcher:02}. Under $\widetilde H_k(B/A)\cong H_k(B,A)$, the splitting defines
\[
  J_k(x)=\left(\bar\pi_*x,\operatorname{susp}^{-1}(q_*x)\right).
\]

It remains to identify the endomorphism. On each interval component of $h^{-1}((-1,1))$, the horizontal direction of $F(\sigma(z))$ is a continuous circle-valued function. It lifts through $\R\to S^1$ to a continuous angle $\beta(z)$ because the component is an interval. By equivariance,
\[
 F(\Phi(z,u))=\Phi(h(z),R_{\beta(z)}u),
\]
where $R_\alpha$ is rotation by $\alpha$. Define the lift
\begin{equation}
  \label{eq:sphere-height-lift}
  F_0(\Phi(z,u))=\Phi(h(z),u)\quad\forall(z,u)\in I\times S^1.
\end{equation}
Replacing $\beta(z)$ by $s\beta(z)$ for $0\leq s\leq1$ joins $F_0$ to $F$. At every remaining height the image is the pole at $h(z)$. The homotopy is continuous there because its horizontal component has norm $\rho(h(z))$, independently of $u$ and $s$. All maps in the homotopy cover $h$, so Lemma~\ref{lem:preimage-pair} shows that the two pointed maps are homotopic.

For the lift in~\eqref{eq:sphere-height-lift}, let $P=(F_0)_{(N,L)}$. Under the identification $W/M\cong(B/K)\wedge S^1$,
\[
  \bar\pi\circ P=\psi_A\circ\bar\pi\quad\text{and}\quad q\circ P=(\psi_K\wedge\mathrm{id}_{S^1})\circ q.
\]
Taking homology and applying the suspension isomorphism identifies the two summands in the conjugacy.
\end{proof}

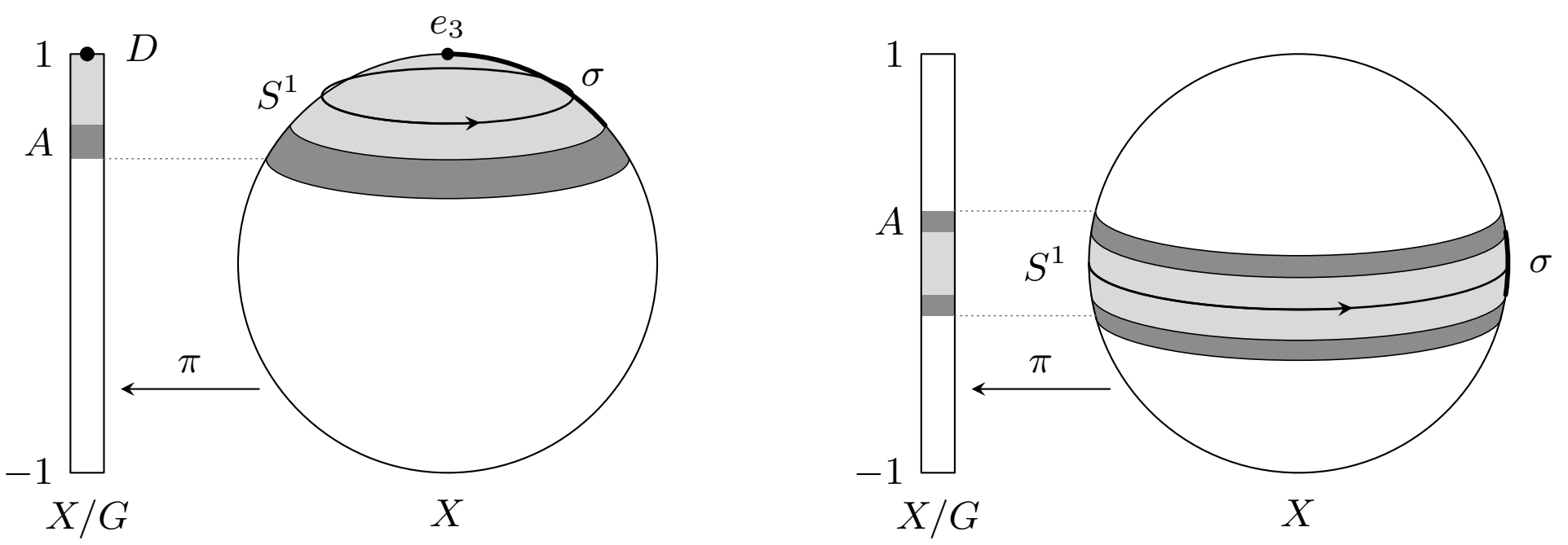
\begin{figure}[!htbp]
\centering
\def\rpecc{0.22}
\newcommand{\rpzone}[3]{%
  \pgfmathsetmacro{\angA}{asin(#1)}\pgfmathsetmacro{\angB}{asin(#2)}%
  \pgfmathsetmacro{\radA}{cos(\angA)}\pgfmathsetmacro{\radB}{cos(\angB)}%
  \fill[#3] (\angA:1) arc[start angle=\angA,end angle=\angB,radius=1]
    arc[start angle=360,end angle=180,x radius=\radB,y radius=\rpecc*\radB]
    arc[start angle=180-\angB,end angle=180-\angA,radius=1]
    arc[start angle=180,end angle=360,x radius=\radA,y radius=\rpecc*\radA] -- cycle;}
\newcommand{\rptopzone}[2]{%
  \pgfmathsetmacro{\angA}{asin(#1)}\pgfmathsetmacro{\radA}{cos(\angA)}%
  \fill[#2] (\angA:1) arc[start angle=\angA,end angle=180-\angA,radius=1]
    arc[start angle=180,end angle=360,x radius=\radA,y radius=\rpecc*\radA] -- cycle;}
\newcommand{\rplatitude}[1]{%
  \pgfmathsetmacro{\radA}{sqrt(1-(#1)*(#1))}%
  \draw[line width=.3pt] (-\radA,#1) arc[start angle=180,end angle=360,x radius=\radA,y radius=\rpecc*\radA];}
\newcommand{\rpframe}{%
  \draw[line width=.4pt] (0,0) circle[radius=1];
  \draw[line width=.4pt] (-1.80,-1) rectangle (-1.64,1);
  \node[left] at (-1.80,1) {$1$};
  \node[left] at (-1.80,-1) {$-1$};
  \node[below] at (-1.72,-1.05) {$X/G$};
  \node[below] at (0,-1.05) {$X$};
  \draw[->,line width=.4pt] (-0.9,-0.6) -- node[above] {$\pi$} (-1.56,-0.6);}
\tikzset{rppairs/.style={scale=1.7,line cap=round,line join=round,>=stealth,every node/.style={font=\footnotesize}}}
\subfloat[At $e_3$]{%
\begin{tikzpicture}[rppairs]
  \rptopzone{0.5}{black!15}
  \rpzone{0.5}{0.66}{black!45}
  \rplatitude{0.5}
  \rplatitude{0.66}
  \fill[black!15] (-1.80,0.5) rectangle (-1.64,1);
  \fill[black!45] (-1.80,0.5) rectangle (-1.64,0.66);
  \rpframe
  \pgfmathsetmacro{\radO}{sqrt(1-0.8*0.8)}
  \draw[line width=.5pt] (0,0.8) ellipse[x radius=\radO,y radius=\rpecc*\radO];
  \draw[line width=.5pt,->] (-\radO,0.8) arc[start angle=180,end angle=285,x radius=\radO,y radius=\rpecc*\radO];
  \node[left] at (-\radO-0.02,0.82) {$S^1$};
  \draw[line width=1.1pt] (41.3:1) arc[start angle=41.3,end angle=90,radius=1];
  \node[right] at (58:1.05) {$\sigma$};
  \fill (0,1) circle[radius=.03];
  \node[above] at (0,1) {$e_3$};
  \fill (-1.72,1) circle[radius=.035];
  \node[right] at (-1.62,1.03) {$D$};
  \node[left] at (-1.80,0.58) {$A$};
  \draw[densely dotted,black!55,line width=.3pt] (-1.64,0.5) -- (-0.866,0.5);
\end{tikzpicture}}\hfil
\subfloat[At the equator]{%
\begin{tikzpicture}[rppairs]
  \rpzone{-0.25}{0.25}{black!15}
  \rpzone{0.15}{0.25}{black!45}
  \rpzone{-0.25}{-0.15}{black!45}
  \rplatitude{0.25}
  \rplatitude{0.15}
  \rplatitude{-0.15}
  \rplatitude{-0.25}
  \fill[black!15] (-1.80,-0.25) rectangle (-1.64,0.25);
  \fill[black!45] (-1.80,0.15) rectangle (-1.64,0.25);
  \fill[black!45] (-1.80,-0.25) rectangle (-1.64,-0.15);
  \rpframe
  \draw[line width=.5pt] (-1,0) arc[start angle=180,end angle=360,x radius=1,y radius=\rpecc];
  \draw[line width=.5pt,->] (-1,0) arc[start angle=180,end angle=285,x radius=1,y radius=\rpecc];
  \node[left] at (-1.02,0) {$S^1$};
  \draw[line width=1.1pt] (-8.63:1) arc[start angle=-8.63,end angle=8.63,radius=1];
  \node[right] at (0:1.02) {$\sigma$};
  \node[left] at (-1.80,0.2) {$A$};
  \draw[densely dotted,black!55,line width=.3pt] (-1.64,0.25) -- (-0.968,0.25);
  \draw[densely dotted,black!55,line width=.3pt] (-1.64,-0.25) -- (-0.968,-0.25);
\end{tikzpicture}}
\caption{Schematic index pairs $(N,L)$ on $X=S^2$ at the north pole $e_3$ and the equator.}
\label{fig:reconstruction-pairs}
\end{figure}

Figure~\ref{fig:reconstruction-pairs} shows the calculation at the repelling north pole and equator of Section~\ref{sec:ex:so2}. At the pole $H_*(B,A)=0$, and the interval from $A$ to $D=\{1\}$ represents the generator of $H_1(B,K)$. Revolving it about the vertical axis gives the cap class in degree $2$. At the equator $D=\emptyset$ and $H_1(B,A)=\F$. The meridian arc and the band it sweeps out contribute in degrees $1$ and $2$. At each attracting latitude the trapping interval misses the poles and has empty exit, so its lift is an annulus with classes in degrees $0$ and $1$. For the reference flow the index map is the identity on each of these groups by Appendix~\ref{supp:sphere-pairs}.

For the vector field $V(p)=a(z)(e_3-zp)+\omega e_3\times p$, let $h_\tau$ be the time-$\tau$ map of $\dot z=a(z)(1-z^2)$. The meridional field commutes with the rotation at angular speed $\omega$, so the time-$\tau$ map sends $\Phi(z,u)$ to $\Phi(h_\tau(z),R_{\omega\tau}u)$. In Section~\ref{sec:ex:so2}, $a(z)=z(0.4^2-z^2)$ and $\omega=0.6$, so the time-ten map rotates by $6$ radians. The learned height map $f^{10}$ fixes the endpoints. Hence~\eqref{eq:sphere-height-lift} with $h=f^{10}$ defines a continuous equivariant lift, and Theorem~\ref{thm:sphere-reconstruction} computes the same index map on the lifted pair for every equivariant lift of $f^{10}$.

\subsection{The \texorpdfstring{$\mathrm{O}(3)$}{O(3)} action on \texorpdfstring{$\R^4$}{R4}}
\label{app:o3-instance}

Let $\mathrm O(3)$ act by $(a,v)\mapsto(a,Qv)$ for $Q\in\mathrm O(3)$. The orbit coordinates are $\pi(a,v)=(a,r)$ with $r=\|v\|$. The orbit at $(a,r)$ is parametrized by $(a,ru)$ for $u\in S^2$ and collapses to a point when $r=0$. Let $H$ denote the half-plane $r\geq0$ for the learned map, or the rectangle $[-3/2,3/2]\times[0,3/2]$ for the physical flow. The corresponding phase space is $X=\pi^{-1}(H)$.

\filbreak
\begin{thm}
\label{thm:orthogonal-reconstruction}
Let $h=(h_1,h_2)\colon H\to H$ be continuous and $(B,A)$ an index pair for $h$ with $B$ nonempty. Suppose $h$ sends the axis $H\cap\{r=0\}$ into itself and $B$ has a finite CW structure in which $A$ and $D=B\cap\{r=0\}$ are subcomplexes. Define $K=A\cup D$, $N=\pi^{-1}(B)$ and $L=\pi^{-1}(A)$. For $\varepsilon\in\{1,-1\}$ the formula
\[
 F_\varepsilon(a,ru)=\left(h_1(a,r),\,\varepsilon h_2(a,r)u\right),\quad (a,r)\in H,\quad u\in S^2,
\]
defines a continuous equivariant lift $F_\varepsilon\colon X\to X$ with index pair $(N,L)$. In every degree $k\geq0$, the index map of $F_\varepsilon$ is conjugate to the direct sum on the right, with the same isomorphism $J_k$ for both signs.
\begin{align*}
  J_k\colon \widetilde H_k(N/L)&\xrightarrow{\cong} H_k(B,A)\oplus H_{k-2}(B,K),\\
  J_k\,((F_\varepsilon)_{(N,L)})_{*,k}\,J_k^{-1}&=(\psi_A)_{*,k}\oplus\varepsilon(\psi_K)_{*,k-2}.
\end{align*}
\end{thm}

\begin{proof}
For $r>0$ the formula is continuous since $u=v/r$. The vector coordinate has norm $h_2(a,r)$, which tends to zero as $(a,r)$ approaches any point of the axis. This proves continuity there. The formula commutes with every $Q\in\mathrm O(3)$ and satisfies $\pi F_\varepsilon=h\pi$. The preimages $N$ and $L$ are compact as the images of $B\times S^2$ and $A\times S^2$. Lemma~\ref{lem:preimage-pair} applied to the constant homotopy shows that $(N,L)$ is an index pair. As in the sphere case, adjoining $D$ to $A$ preserves the two exit inclusions, so $\psi_K$ is defined.

Choose the section $(a,r)\mapsto(a,re_1)$, and let $M$ be its pointed image in $W=N/L$. Give the sphere of directions the vertex $e_1$ and one open $2$-cell. The product cell construction in the proof of Theorem~\ref{thm:sphere-reconstruction} makes $M\cong B/A$ a subcomplex. Collapsing $M$ now identifies $W/M$ with $(B/K)\wedge S^2=\Sigma^2(B/K)$. Projection onto $B/A$ and collapse of $M$ split the homology of $W$ as before. For $\varepsilon=1$ the lift preserves the section and its pointed map on $W/M$ is $\psi_K\wedge\mathrm{id}_{S^2}$.

For $\varepsilon=-1$, let $R_\theta$ denote rotation about $e_3$ and define
\[
  F^{(s)}(a,ru)=\left(h_1(a,r),h_2(a,r)R_{s\pi}(-u)\right)\quad\forall s\in[0,1].
\]
Then $F^{(0)}=F_{-1}$ and $\pi\circ F^{(s)}=h\circ\pi$ for every $s$. The vector coordinate has norm $h_2(a,r)$, so the homotopy is continuous at $r=0$. Lemma~\ref{lem:preimage-pair} therefore shows that the endpoint pointed maps are homotopic. Since $R_\pi(-u)=T(u)$ for the reflection $T(u_1,u_2,u_3)=(u_1,u_2,-u_3)$ and $T(e_1)=e_1$, the final lift acts by $\psi_A$ on the section and by $\psi_K\wedge T$ on the collapsed space. On $(B/K)\wedge S^2$, the map $\mathrm{id}\wedge T$ multiplies each reduced cellular chain by $\deg T=-1$, since $T$ reverses the orientation of the open $2$-cell. Composing with $\psi_K\wedge\mathrm{id}_{S^2}$ proves the formula for the second summand.
\end{proof}

For the physical flow of Section~\ref{sec:ex:o3}, $\dot v$ is a scalar multiple of $v$, so the time-$\log3$ map has $\varepsilon=1$. The quotient pair at the origin in Appendix~\ref{supp:o3-pairs} is a half-disk with a semicircular exit collar. Subdividing along the collar and the diameter makes $A$ and $D$ subcomplexes. The saddle and sink pairs there are rectangles. Composing the time map with $(a,v)\mapsto(a,-v)$ leaves the quotient dynamics unchanged and reverses the sign on the shifted summand.

\paragraph{Learned lift and local pairs.}
The learned map has the form $f(s,r)=(s g_1(s,r),r g_2(s,r))$, where $s\in\R$, $r\geq0$, and the learned gains $g_1,g_2$ are positive. Its specified lift on $\R\times\R^3$ is
\[
  \widetilde F(s,w)=(s g_1(s,\|w\|),g_2(s,\|w\|)w).
\]
Thus Theorem~\ref{thm:orthogonal-reconstruction} applies with $h=f$ and $\varepsilon=1$. The lift is defined by the learned gains independently of the decoder.

For the selected fit, interval bounds verify expansion in both quotient coordinates at the origin, expansion in $s$ and contraction on $s=0$ at the saddle, and contraction in $r$ and along $r=0$ at each sink. On each of the verified boxes let $g$ be a diagonal affine map centered at the fixed point, with positive multipliers greater than one in the expanding directions and less than one in the contracting directions. The same inequalities hold along $(1-t)f+tg$, $0\leq t\leq1$. At the origin every nonzero coordinate expands. At the saddle, expansion forces every full trajectory in the box onto $s=0$, where contraction forces it to be the fixed point. At a sink, contraction in $r$ forces every full trajectory onto $r=0$, and contraction there forces the fixed point. These boxes therefore remain isolating, as do their preimages under $\pi$, and continuation preserves the Conley indices of the specified lifts~\citep[Theorem~4.9]{franks:richeson:00}.

For the diagonal normal map, take a small rectangle about the fixed point. In each expanding direction normalize the coordinate to $x\in[-1,1]$ (or $x\in[0,1]$ for the radius at the origin), and let $\lambda>1$ be its multiplier. The discrete exit set is the union of the collars $|x|\geq1/\lambda$ in these coordinates. The same exit inclusions hold for intermediate multipliers $\lambda^t$, $0\leq t\leq1$, with contracting multipliers also raised to the power $t$. Thus the index maps of the normal lifts preserving direction are identities. In each expanding coordinate, rescale the central interval linearly and collapse the collar to its outer endpoint. The straight-line homotopy to this coordinate map preserves the collar, its outer faces, and the bottom axis when present. The resulting pairs with exit faces therefore compute both $H_*(B,A)$ and $H_*(B,K)$.

At the origin the exit faces are the two sides and the top. They form an arc, and adjoining the bottom axis gives the whole boundary, so only $H_2(B,K)$ contributes and its class shifts to degree $4$. At the saddle the two exit sides give $H_1(B,A)=\F$, and $D=\emptyset$, so the two classes are in degrees $1$ and $3$. At each sink $A=\emptyset$ and $D$ is the bottom edge, so only $H_0(B,A)=\F$ contributes. The relative-chain computation uses $\mathbb F_5$. Reversing the vector direction replaces the identity by $-1$ in degree $4$ at the origin and degree $3$ at the saddle. Figure~\ref{fig:o3}(d) attaches these local indices to the grouped sampled graph. Its connecting arrows and the correspondence with the physical flow are separate numerical comparisons.

\subsection{Reduced attitude control on \texorpdfstring{$S^2$}{S2}}
\label{app:attitude-instance}

On the invariant subsystem $\{\eta=0\}$ of Section~\ref{sec:ex:attitude}, the thrust axis follows $V_a(q)=a(z)(e_3-zq)+\omega e_3\times q$, with $z=q\cdot e_3$ and $\omega=\pi/\log2$.

\begin{cor}
\label{cor:attitude-axis-index}
For the time-one map of the thrust axis, the index maps on the pairs of Appendix~\ref{supp:attitude-pairs} are identities. With nominal gain $a=1$, the nonzero groups are $\F$ in degree $0$ at $e_3$ and $\F$ in degree $2$ at $-e_3$. With coning gain $a(z)=1/4-z$, each pole has $\F$ in degree $2$, and the attracting latitude $C=\{z=1/4\}$ has $\F$ in each of degrees $0$ and $1$. All other groups vanish.
\end{cor}

\begin{proof}
The height satisfies $\dot z=a(z)(1-z^2)$. Between the poles it increases for $a=1$ and moves toward $1/4$ for $a(z)=1/4-z$. Apply Theorem~\ref{thm:sphere-reconstruction} to the height projections of the stated pairs. At an attracting pole only $H_0(B,A)=\F$ contributes. At the attracting latitude $A=D=\emptyset$, so the two copies of $H_0(B)=\F$ occur in degrees $0$ and $1$. At a repelling pole only $H_1(B,K)=\F$ contributes, in degree $2$. The exit inclusions hold for $0\leq s\leq1$, so Lemma~\ref{lem:pointed-homotopy} gives identity index maps.
\end{proof}

\section{Analytic and computational details}

This appendix follows the examples from the equations to the computed graphs. The analytical panels provide a reference, the learned panels show the numerical model, and the reconstructed panels describe indices in a specified lifted state space. These are different objects. In particular, agreement of sampled graphs is not a uniform error bound for the original dynamics.

\subsection{Sphere example}
\label{supp:sphere-walkthrough}
\begin{figure}[!htb]
\centering
\begin{minipage}[b]{.49\linewidth}\centering
\subfloat[Sphere trajectories]{\includegraphics[height=2.05in]{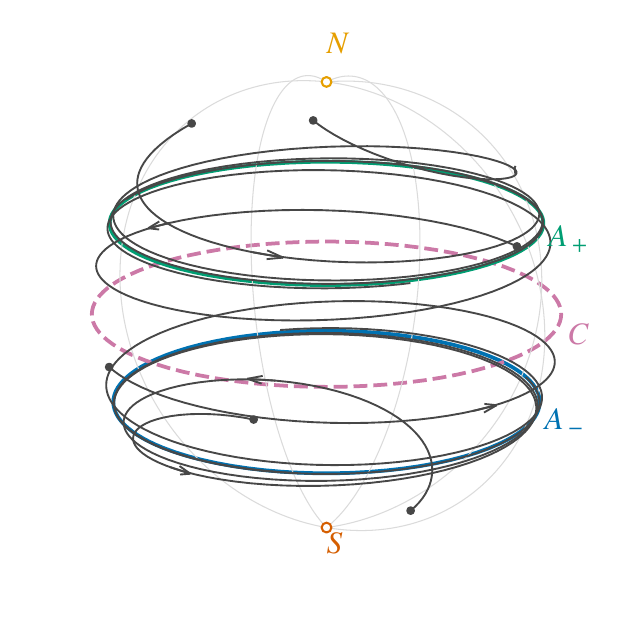}}
\end{minipage}\hfill
\begin{minipage}[b]{.49\linewidth}\centering
\subfloat[Physical height map]{\includegraphics[height=2.05in]{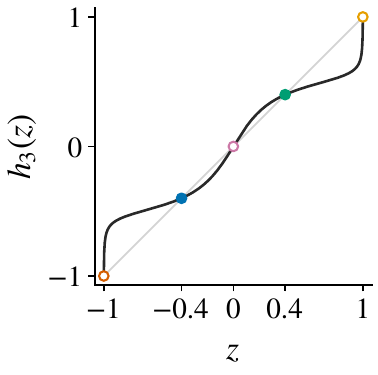}}
\end{minipage}
\caption{Sphere flow and its height quotient.}
\label{fig:appendix-sphere-height}
\end{figure}

The height coordinate identifies $S^2/G$ with $[-1,1]$, where the quotient flow satisfies $\dot z=z(0.4^2-z^2)(1-z^2)$. Figure~\ref{fig:appendix-sphere-height} shows trajectories on the sphere and the time-three map of this equation. The model in Section~\ref{sec:ex:so2} is trained on time-one pairs. The graph comparison uses its tenth iterate $f^{10}$ and the time-ten quotient map $\olF$.

\begin{figure}[!htbp]
\centering
\noindent\begin{minipage}[b]{.485\linewidth}\centering
Analytical\par\smallskip
\includegraphics[height=1.5in]{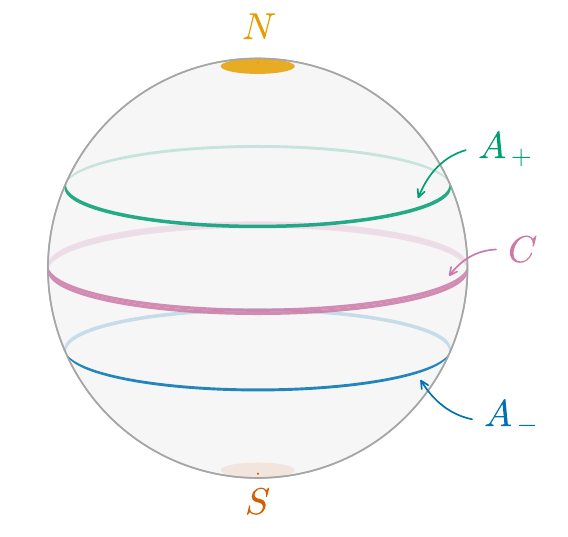}
\end{minipage}\hfill
\begin{minipage}[b]{.485\linewidth}\centering
Learned\par\smallskip
\includegraphics[height=1.5in]{example1/learned-morse-preimages.pdf}
\end{minipage}\par\vspace{2pt}
\panelpair{example1/reference-morse-interval.pdf}{example1/learned-morse-interval.pdf}
\panelpair{example1/reference-quotient.pdf}{example1/learned-quotient-appendix.pdf}
\panelpair{example1/reference-lift.pdf}{example1/learned-lift-appendix.pdf}
\caption{Sphere example. Analytical (left) and learned (right).
Rows show Morse-set preimages on the sphere, quotient intervals,
quotient Conley--Morse graphs, and reconstructed sphere indices.
Learned bands are the saved Morse-cell preimages, colored by group.
The learned local indices refer to the specified sphere lift.}
\label{fig:appendix-sphere-comparison}
\end{figure}

The encoder preimages locate the learned Morse sets on the sphere. Figure~\ref{fig:appendix-sphere-comparison} compares them with the preimages of the Morse sets computed for $\olF$. Both interval computations use $2^7$ cells. The learned computation has six Morse sets. Two adjacent Morse sets representing the learned counterpart of $C$ are combined, and CMGDB recomputes the index on the merged index pair. The resulting five-node graph has the same directed connections and index annotations as the analytical graph. The reconstructed sphere indices agree with the analytical indices derived from the explicit index pairs in Appendix~\ref{supp:sphere-pairs}.

\paragraph{Training.}
The encoder, latent map, and orbit decoder are trained jointly using $4096$ training starts and $20000$ updates on time-one pairs. With the encoder and decoder fixed, the latent map is then trained for $4000$ further updates on the sum of the mean squared ten-step residual, the mean of the $32$ largest squared ten-step residuals in each batch, and $0.1$ times the mean squared one-step residual. The checkpoint with the smallest validation maximum of the ten-step residual, reached after $3400$ of these updates, is kept. The graph computation uses the tenth iterate of this map. Its networks use width $32$, two layers, batches of $256$, and a learning rate decreasing from $10^{-3}$ to $2\cdot10^{-5}$. The cycle and spread weights are $0.1$ and $0.01$.

\paragraph{Index pairs.}
\label{supp:sphere-pairs}

By continuity of $h_{10}$ at its fixed points $1$ and $0$, there are $1/2<c<1$ and $0<c'<1/4$ with $h_{10}(c)\geq1/2$ and $h_{10}(c')\leq1/4$. Consider the following pairs in height coordinates.
\[
  \begin{array}{c|cc}
    & B & A \\ \hline
    e_3 & [1/2,1] & [1/2,c] \\
    C & [-1/4,1/4] & [-1/4,-c']\cup[c',1/4] \\
    A_+ & [1/4,3/4] & \emptyset
  \end{array}
\]
Since $h_s(c)$ decreases and $h_s(c')$ increases with $s$, the lifted pairs $(N,L)=(\pi^{-1}(B),\pi^{-1}(A))$ satisfy both exit inclusions for $0\leq s\leq10$ and isolate $e_3$, $C$ and $A_+$, respectively. They are homotopy equivalent to $(D^2,S^1)$, $S^1\times([-1,1],\setof{\pm1})$ and $(S^1,\emptyset)$, so for each pair $(N,L)$ the index map is the identity on $H_k(N,L;\F)$, which is isomorphic to $\F$ in degree $2$ at $e_3$, degrees $1$ and $2$ at $C$, and degrees $0$ and $1$ at $A_+$, and vanishes in all other degrees. By the symmetry $z\mapsto-z$, the same conclusions hold at $-e_3$ and $A_-$. The formula of Theorem~\ref{thm:sphere-reconstruction} returns the same indices (Figure~\ref{fig:sphere-reconstruction}).

\begin{figure}[H]
\centering
\includegraphics[height=1.5in]{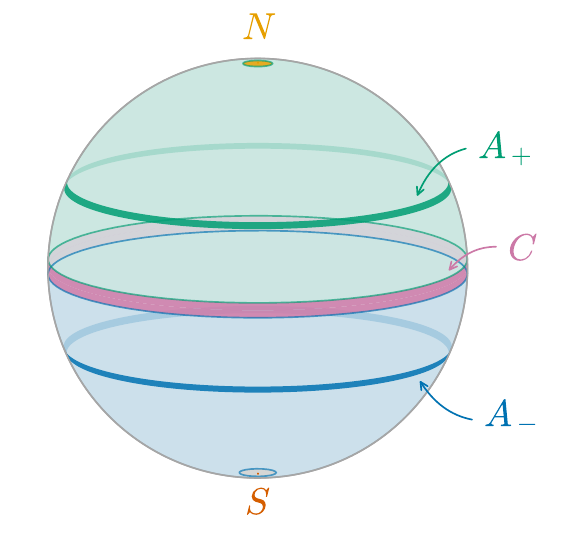}
\caption{Encoder preimages of Morse sets (dark bands) and attracting blocks (pale regions).}
\label{fig:appendix-sphere-blocks}
\end{figure}

\paragraph{Attracting blocks.}
\label{supp:sphere-computation}
On the $2^7$-cell grid, the maximum sampled residuals are $0.00726$ and $0.00738$, below the estimated tolerances $0.04137$ and $0.05667$, respectively. The comparison uses the ten-step learned map and the corresponding blocks associated with $A_-$ and $A_+$, consisting of $57$ and $67$ cells from which only the respective attracting Morse set is reachable. Figure~\ref{fig:appendix-sphere-blocks} shows their encoder preimages. The residual maxima and tolerances are sampled estimates.

\FloatBarrier

\subsection{Orthogonal example}
\label{supp:o3-walkthrough}
Write a state as $(a,v)\in\R\times\R^3$ and let $r=\|v\|$. On $X=[-3/2,3/2]\times B^3_{3/2}$, the quotient map $q:X\to[-3/2,3/2]\times[0,3/2]$ is $q(a,v)=(a,r)$. The quotient equations are $\dot a=a(1-a^2)$ and $\dot r=(1-2a^2-r^2)r$. The graph calculation uses their time-$\log3$ map. Figure~\ref{fig:o3}(a) shows the physical quotient flow; each point with $r>0$ represents an entire spherical orbit.

In Figure~\ref{fig:appendix-o3-comparison}, the colored cells are numerical quotient regions. Finite-graph basin labels describe reachability in the computed graph, not a certified basin of the original differential equation. The bottom row concerns the specified equivariant lift. It should not be read as an index computed directly from a four-dimensional plot.
\usetikzlibrary{shapes.geometric}
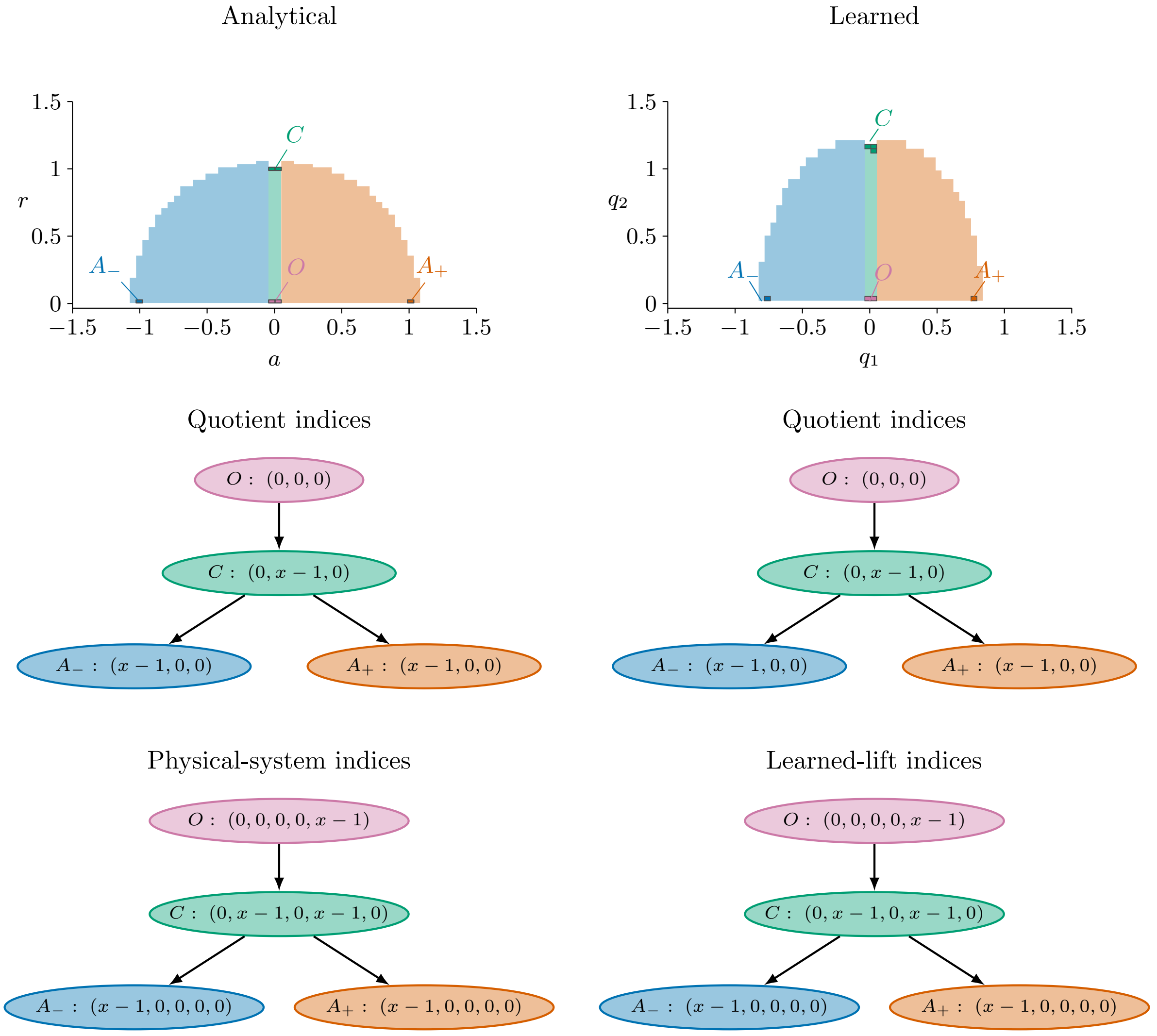
\begin{figure}[ht]
\centering
\begingroup
\definecolor{oA}{rgb}{0,0.447059,0.698039}
\definecolor{oB}{rgb}{0.835294,0.368627,0}
\definecolor{oC}{rgb}{0,0.619608,0.45098}
\definecolor{oO}{rgb}{0.8,0.47451,0.654902}
\definecolor{oAf}{rgb}{0.6,0.780392,0.878431}
\definecolor{oBf}{rgb}{0.933333,0.74902,0.6}
\definecolor{oCf}{rgb}{0.6,0.847059,0.780392}
\definecolor{oOf}{rgb}{0.921569,0.788235,0.862745}
\newcommand{\oGraph}[1]{%
\begin{tikzpicture}[x=1bp,y=1bp,>=latex,
 every node/.style={font=\fontsize{7}{8}\selectfont,inner sep=2bp},
 nd/.style={ellipse,draw,line width=.8bp,minimum height=17bp,align=center},
 every path/.style={line width=.8bp}]
\path[use as bounding box] (0,-16) rectangle (210,93);
\ifnum#1=0
\node[nd,draw=oO,fill=oOf,minimum width=65bp] (O) at (105,78) {};
\node[inner sep=0] at (O) {$O:\;(0,0,0)$};
\node[nd,draw=oC,fill=oCf,minimum width=90bp] (C) at (105,42) {};
\node[inner sep=0] at (C) {$C:\;(0,x-1,0)$};
\node[nd,draw=oA,fill=oAf,minimum width=90bp] (A) at (49,6) {};
\node[inner sep=0] at (A) {$A_-:\;(x-1,0,0)$};
\node[nd,draw=oB,fill=oBf,minimum width=90bp] (B) at (161,6) {};
\node[inner sep=0] at (B) {$A_+:\;(x-1,0,0)$};
\else
\node[nd,draw=oO,fill=oOf,minimum width=100bp] (O) at (105,78) {};
\node[inner sep=0] at (O) {$O:\;(0,0,0,0,x-1)$};
\node[nd,draw=oC,fill=oCf,minimum width=100bp] (C) at (105,42) {};
\node[inner sep=0] at (C) {$C:\;(0,x-1,0,x-1,0)$};
\node[nd,draw=oA,fill=oAf,minimum width=100bp] (A) at (49,6) {};
\node[inner sep=0] at (A) {$A_-:\;(x-1,0,0,0,0)$};
\node[nd,draw=oB,fill=oBf,minimum width=100bp] (B) at (161,6) {};
\node[inner sep=0] at (B) {$A_+:\;(x-1,0,0,0,0)$};
\fi
\draw[->] (O)--(C);\draw[->] (C)--(A);\draw[->] (C)--(B);
\end{tikzpicture}}
\makebox[.485\linewidth]{\small Analytical}\hfill\makebox[.485\linewidth]{\small Learned}\par
\makebox[.485\linewidth]{\begin{tikzpicture}[x=1bp,y=1bp]
\path[use as bounding box] (0,0) rectangle (210,140);
\node[anchor=south west,inner sep=0] at (4.108460,15.710874) {\includegraphics[width=185.533981bp,height=108.432501bp]{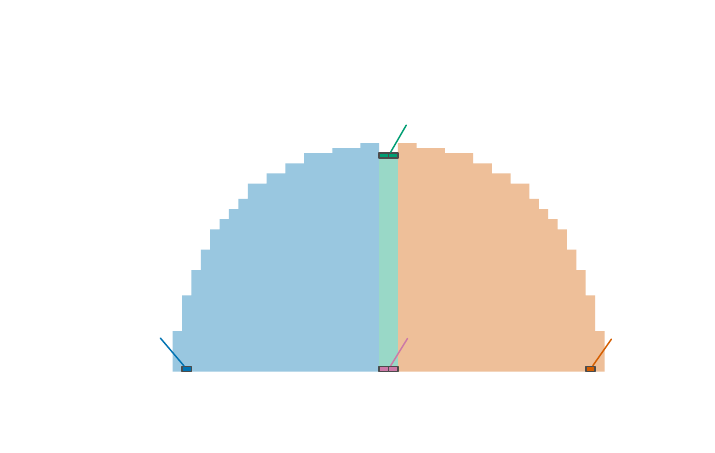}};
\node[inner sep=0,text=oC,font=\fontsize{9}{10}\selectfont] at (113.048821,100.308213) {$C$};
\node[inner sep=0,text=oO,font=\fontsize{9}{10}\selectfont] at (113.481553,49.090096) {$O$};
\node[inner sep=0,text=oA,font=\fontsize{9}{10}\selectfont] at (39.754785,49.040127) {$A_-$};
\node[inner sep=0,text=oB,font=\fontsize{9}{10}\selectfont] at (166.166713,49.040127) {$A_+$};
\draw[line width=.4bp] (27,33.18)--(183,33.18);\draw[line width=.4bp] (27,33.18)--(27,113);
\draw[line width=.4bp] (27.0,33.18)--(27.0,31.18);\node[font=\fontsize{9}{10}\selectfont,anchor=north,inner sep=1bp] at (27.0,30) {$-1.5$};
\draw[line width=.4bp] (53,33.18)--(53,31.18);\node[font=\fontsize{9}{10}\selectfont,anchor=north,inner sep=1bp] at (53,30) {$-1$};
\draw[line width=.4bp] (79.0,33.18)--(79.0,31.18);\node[font=\fontsize{9}{10}\selectfont,anchor=north,inner sep=1bp] at (79.0,30) {$-0.5$};
\draw[line width=.4bp] (105,33.18)--(105,31.18);\node[font=\fontsize{9}{10}\selectfont,anchor=north,inner sep=1bp] at (105,30) {$0$};
\draw[line width=.4bp] (131.0,33.18)--(131.0,31.18);\node[font=\fontsize{9}{10}\selectfont,anchor=north,inner sep=1bp] at (131.0,30) {$0.5$};
\draw[line width=.4bp] (157,33.18)--(157,31.18);\node[font=\fontsize{9}{10}\selectfont,anchor=north,inner sep=1bp] at (157,30) {$1$};
\draw[line width=.4bp] (183.0,33.18)--(183.0,31.18);\node[font=\fontsize{9}{10}\selectfont,anchor=north,inner sep=1bp] at (183.0,30) {$1.5$};
\draw[line width=.4bp] (27,35)--(25,35);\node[font=\fontsize{9}{10}\selectfont,anchor=east,inner sep=1bp] at (24,35) {$0$};
\draw[line width=.4bp] (27,61.0)--(25,61.0);\node[font=\fontsize{9}{10}\selectfont,anchor=east,inner sep=1bp] at (24,61.0) {$0.5$};
\draw[line width=.4bp] (27,87)--(25,87);\node[font=\fontsize{9}{10}\selectfont,anchor=east,inner sep=1bp] at (24,87) {$1$};
\draw[line width=.4bp] (27,113.0)--(25,113.0);\node[font=\fontsize{9}{10}\selectfont,anchor=east,inner sep=1bp] at (24,113.0) {$1.5$};
\node[font=\fontsize{9}{10}\selectfont] at (105,13) {$a$};\node[font=\fontsize{9}{10}\selectfont] at (8,74) {$r$};
\end{tikzpicture}
}\hfill
\makebox[.485\linewidth]{\begin{tikzpicture}[x=1bp,y=1bp]
\path[use as bounding box] (0,0) rectangle (210,140);
\node[anchor=south west,inner sep=0] at (34.961349,10.210256) {\includegraphics[width=129.276883bp,height=99.312410bp]{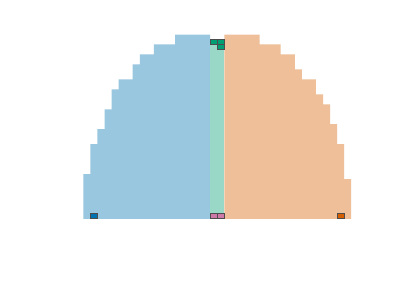}};
\node[inner sep=0,text=oC,font=\fontsize{9}{10}\selectfont] at (110.370296,106.717746) {$C$};
\node[inner sep=0,text=oO,font=\fontsize{9}{10}\selectfont] at (110.370296,47.113198) {$O$};
\node[inner sep=0,text=oA,font=\fontsize{9}{10}\selectfont] at (56.506762,47.113198) {$A_-$};
\node[inner sep=0,text=oB,font=\fontsize{9}{10}\selectfont] at (151.441240,47.113198) {$A_+$};
\draw[oC,line width=.4bp] (108.350413,102.510367) -- (104.983943,97.601757);
\draw[oO,line width=.4bp] (108.350413,43.607048) -- (104.983943,35.893518);
\draw[oA,line width=.4bp] (58.526645,43.607048) -- (63.239704,35.893518);
\draw[oB,line width=.4bp] (149.421358,43.607048) -- (144.035005,35.893518);
\draw[line width=.4bp] (27,33.18)--(183,33.18);\draw[line width=.4bp] (27,33.18)--(27,113);
\draw[line width=.4bp] (27.0,33.18)--(27.0,31.18);\node[font=\fontsize{9}{10}\selectfont,anchor=north,inner sep=1bp] at (27.0,30) {$-1.5$};
\draw[line width=.4bp] (53,33.18)--(53,31.18);\node[font=\fontsize{9}{10}\selectfont,anchor=north,inner sep=1bp] at (53,30) {$-1$};
\draw[line width=.4bp] (79.0,33.18)--(79.0,31.18);\node[font=\fontsize{9}{10}\selectfont,anchor=north,inner sep=1bp] at (79.0,30) {$-0.5$};
\draw[line width=.4bp] (105,33.18)--(105,31.18);\node[font=\fontsize{9}{10}\selectfont,anchor=north,inner sep=1bp] at (105,30) {$0$};
\draw[line width=.4bp] (131.0,33.18)--(131.0,31.18);\node[font=\fontsize{9}{10}\selectfont,anchor=north,inner sep=1bp] at (131.0,30) {$0.5$};
\draw[line width=.4bp] (157,33.18)--(157,31.18);\node[font=\fontsize{9}{10}\selectfont,anchor=north,inner sep=1bp] at (157,30) {$1$};
\draw[line width=.4bp] (183.0,33.18)--(183.0,31.18);\node[font=\fontsize{9}{10}\selectfont,anchor=north,inner sep=1bp] at (183.0,30) {$1.5$};
\draw[line width=.4bp] (27,35)--(25,35);\node[font=\fontsize{9}{10}\selectfont,anchor=east,inner sep=1bp] at (24,35) {$0$};
\draw[line width=.4bp] (27,61.0)--(25,61.0);\node[font=\fontsize{9}{10}\selectfont,anchor=east,inner sep=1bp] at (24,61.0) {$0.5$};
\draw[line width=.4bp] (27,87)--(25,87);\node[font=\fontsize{9}{10}\selectfont,anchor=east,inner sep=1bp] at (24,87) {$1$};
\draw[line width=.4bp] (27,113.0)--(25,113.0);\node[font=\fontsize{9}{10}\selectfont,anchor=east,inner sep=1bp] at (24,113.0) {$1.5$};
\node[font=\fontsize{9}{10}\selectfont] at (105,13) {$q_1$};\node[font=\fontsize{9}{10}\selectfont] at (8,74) {$q_2$};
\end{tikzpicture}
}\par
\makebox[.485\linewidth]{\small Quotient indices}\hfill\makebox[.485\linewidth]{\small Quotient indices}\par\smallskip
\makebox[.485\linewidth]{\oGraph{0}}\hfill\makebox[.485\linewidth]{\oGraph{0}}\par\medskip
\makebox[.485\linewidth]{\small Physical-system indices}\hfill\makebox[.485\linewidth]{\small Learned-lift indices}\par\smallskip
\makebox[.485\linewidth]{\oGraph{1}}\hfill\makebox[.485\linewidth]{\oGraph{1}}
\endgroup
\caption{$\mathrm O(3)$ example. Analytical (left) and learned (right). Rows show quotient Morse sets and finite-graph basin labels, grouped quotient graphs with index annotations, and physical-system or learned-lift index annotations. The learned indices refer to Appendix~\ref{app:o3-instance}.}
\label{fig:appendix-o3-comparison}
\end{figure}

\paragraph{Grouping.}
\label{supp:o3-grouping}
On the $2^6\times2^6$ rectangular grid, the six computed Morse sets are combined into four groups: two Morse sets each at $O$ and $C$, and one at each attractor $A_\pm$. Each merged region includes the cells on directed paths between its Morse sets. CMGDB checks the merged index pairs and recomputes their indices on the rectangular domain.

\paragraph{Training.}
The encoder, latent map, and orbit decoder in Section~\ref{sec:ex:o3} are trained jointly. Training uses $16384$ training starts, $4096$ validation starts, $4096$ test starts, and $20000$ updates. The networks have width $32$, two layers, batches of $256$, and a learning rate decreasing from $10^{-3}$ to $2\cdot10^{-5}$. The cycle and spread weights are $0.1$ and $0.01$. Its gain parameterization uses time-$\log3$ pairs and a one-step graph map.

\paragraph{Index pairs.}
\label{supp:o3-pairs}

Both $a^2$ and $a^2+\|v\|^2$ satisfy $\dot y=2y(1-y)$, whose time-$\log3$ map $y\mapsto9y/(1+8y)$ sends $1/8$ to $9/16$ and $1/136$ to $1/16$. Moreover, $\|v\|$ decreases where $|a\mp1|\leq1/8$, and where $|a|\leq1/4$ it increases at $\|v\|=3/4$ and decreases at $\|v\|=5/4$. So the pairs
\begin{align*}
  N&=\setof{a^2+\|v\|^2\leq9/16}, & L&=N\cap\setof{a^2+\|v\|^2\geq1/8} & &\text{at }O,\\
  N&=\setof{|a|\leq1/4,\ 3/4\leq\|v\|\leq5/4}, & L&=N\cap\setof{a^2\geq1/136} & &\text{at }C,\\
  N&=\setof{|a\mp1|\leq1/8,\ \|v\|\leq1/4}, & L&=\emptyset & &\text{at }A_\pm,
\end{align*}
satisfy both inclusions for $0\leq s\leq\log3$ and isolate $O$, $C$ and $A_\pm$. They are homotopy equivalent to $(D^4,S^3)$, $([-1,1],\setof{\pm1})\times S^2$ and a point. For each pair $(N,L)$, the index map is the identity on $H_k(N,L;\F)$, which is isomorphic to $\F$ in degree $4$ at $O$, degrees $1$ and $3$ at $C$, and degree $0$ at $A_\pm$, and vanishes in all other degrees.

\FloatBarrier

\subsection{Attitude example}
\paragraph{Training.}
The encoder receives the 18 coordinates consisting of $p$, $v$, the three columns of $R$, and $\Omega$. Training includes a mean-squared penalty between the encoded state and the known thrust direction $q=Re_3$, whose input coordinates are left unchanged by normalization.

\phantomsection\label{supp:attitude-grouping}
For the coning gain, we merge the additional one-cell Morse set with the Morse set labeled $C$ and all cells on directed paths between them. CMGDB checks the merged index pair and recomputes its interval index as $(x-1,0)$, yielding the graph in Figure~\ref{fig:attitude-computation}(d).

\paragraph{Orbital stability.}
For the coning gain, the orbit $C\times\{0\}$ in $(q,\eta)$ has period $T=2\pi/\omega=2\log2$. On $\eta=0$, the azimuth advances at constant speed $\omega$, while the derivative of $(1/4-z)(1-z^2)$ at $z=1/4$ is $-15/16$. Since $\dot\eta=-2\eta$ is independent of $q$, the variational system along the orbit is block triangular. Its transverse Floquet multipliers are therefore $\exp(-15T/16)=\exp(-15\log2/8)$ and $\exp(-2T)=1/16$, the latter with multiplicity three. The remaining multiplier is $1$, corresponding to phase along the orbit. All transverse multipliers lie inside the unit circle, so the orbit is locally exponentially orbitally stable in $(q,\eta)$.

\paragraph{Index pairs.}
\label{supp:attitude-pairs}

The time-$s$ maps $h_s$ of the height equation are increasing and fix its equilibria. At a repelling pole of height $\varepsilon\in\{-1,1\}$, continuity at $\varepsilon$ gives $1/2<c<1$ with $\varepsilon h_1(\varepsilon c)>1/2$. On the cap $N=\setof{\varepsilon z\geq1/2}$, $\varepsilon z$ decreases away from the pole. Thus $(N,L)$ with $L=\setof{1/2\leq\varepsilon z\leq c}$ satisfies both exit inclusions for $0\leq s\leq1$ and isolates the pole. With empty exit, the cap $N=\setof{z\geq0}$ at the nominal sink $e_3$ and the band $N=\setof{0\leq z\leq1/2}$ at the coning latitude $C$ map into their interiors for every $s>0$. Corollary~\ref{cor:attitude-axis-index} computes the indices of these pairs.

\paragraph{Attracting blocks.}
\label{supp:attitude-computation}

Using the chordal metric on the latent sphere $Z$, we evaluate the residual modulo $G$ of $f^{10}$ over one second on plant trajectories initialized in $\setof{\eta=0}$. For nominal and coning control, the maximum sampled residuals are both $0.015$, below the estimated attracting-block tolerances $0.29$ at the nominal upright attractor and $0.17$ at the grouped coning region, respectively. These sampled comparisons concern the invariant subsystem $\setof{\eta=0}$ and do not provide uniform residual bounds.

\FloatBarrier

\subsection{Reference index maps and timings}
\label{supp:example-pairs}

For each reference time map $F=\varphi_t$, Proposition~\ref{prop:invariant-sets}(ii) lifts the quotient index pair to $(N,L)$. On the pairs of Appendices~\ref{supp:sphere-pairs}, \ref{supp:o3-pairs} and~\ref{supp:attitude-pairs}, the quotient flow satisfies both exit inclusions for $0\leq s\leq t$. Taking preimages gives $\varphi_s(L)\cap N\subset L$ and $\varphi_s(N\setminus L)\subset N$. Lemma~\ref{lem:pointed-homotopy} gives a pointed homotopy from the identity on $N/L$ to $F_{(N,L)}$, so the index map is the identity on $H_*(N,L;\F)$ for these good or empty-exit pairs. The reconstructed graphs use indices of learned lifts in Sections~\ref{sec:ex:so2} and~\ref{sec:ex:o3} and reference thrust-axis indices in Section~\ref{sec:ex:attitude}. In each example the correspondence of learned nodes with physical invariant sets is observed.

\paragraph{Timings.}
\label{supp:timings}

Table~\ref{tab:timing} compares CMGDB times on the quotient grids of $2^7$, $2^5\times2^5$ and $2^7$ cells of Sections~\ref{sec:ex:so2}, \ref{sec:ex:o3} and~\ref{sec:ex:attitude} with those on a grid of a box containing $X$ at the same cell side. CMGDB needs a box as its domain, so the rows of Sections~\ref{sec:ex:so2} and~\ref{sec:ex:attitude} use maps built from samples, of the time-ten map and of the time-one map of the coning gain, each without its rotation about $e_3$, in $[-1,1]^3$ and of their quotient maps in $[-1,1]$, with cubes without samples sent outside $[-1,1]^3$, and the $\mathrm O(3)$ rows use floating-point hulls. For the sphere flow, CMGDB returns in $[-1,1]^3$ the same Morse graph as on the quotient, no index at $A_\pm$ and $C$, and zero index at the poles, which are on faces of $[-1,1]^3$. In $[-3/2,3/2]^4$, it returns $164$ Morse sets, four with the indices of Section~\ref{sec:ex:o3} and $160$ with zero index. For the coning gain, it returns in $[-1,1]$ the Morse sets of $\pm e_3$ and $C$ with the indices of Figure~\ref{fig:attitude-computation}(d) and two one-cell Morse sets with zero index beside $C$. In $[-1,1]^3$, it returns $13$ Morse sets: $C$ without an index, two bands beside $C$ and the cells at $-e_3$ with zero index, the cells at $e_3$ with index $\F$ in degree $2$, and eight small Morse sets with zero index around them. The outputs of the maps built from samples belong to those maps, not to $F$. The plant itself is not computed: a uniform grid with $2^7$ intervals on each of its $12$ coordinates has $2^{84}\approx1.9\cdot10^{25}$ cells.

\begin{table}[!htbp]
  \centering
  \caption{CMGDB timings. In the original space, Cells counts cells of the grid of the bounding box, not only those that meet $X$: all of them on the uniform rows, and those of the final grid on the adaptive rows. The blocks are Sections~\ref{sec:ex:so2}, \ref{sec:ex:o3} and~\ref{sec:ex:attitude}. The 12-dimensional row counts a grid of the plant and is not computed; the attitude rows were measured in a separate run.}
  \label{tab:timing}
  \setlength{\tabcolsep}{5.5pt}
  \begin{tabular}{cccrrrrr}
    Dimension & $G$ & Grid & Cells & Map & Morse graph & Conley index & Total \\
    \hline
    $3$ & \multirow{3}{*}{$\mathrm{SO}(2)$} & uniform & $2^{21}$ & 31\,s & 1.1\,s & 14\,s & 47\,s \\
    $3$ &  & adaptive & 17,464 & 2.9\,s & 60\,ms & 21\,s & 24\,s \\
    \bfseries\boldmath $1$ &  & \bfseries\boldmath uniform & \bfseries\boldmath $2^7$ & \bfseries\boldmath 2.3\,ms & \bfseries\boldmath $<$1\,ms & \bfseries\boldmath $<$1\,ms & \bfseries\boldmath 3.1\,ms \\
    \hline
    $4$ & \multirow{3}{*}{$\mathrm{O}(3)$} & uniform & $2^{23}$ & 16\,s & 4.9\,s & 991\,s & 1,011\,s \\
    $4$ &  & adaptive & 109,536 & 2.3\,s & 580\,ms & 943\,s & 946\,s \\
    \bfseries\boldmath $2$ &  & \bfseries\boldmath uniform & \bfseries\boldmath $2^{10}$ & \bfseries\boldmath $<$1\,ms & \bfseries\boldmath $<$1\,ms & \bfseries\boldmath $<$1\,ms & \bfseries\boldmath 1.8\,ms \\
    \hline
    $12$ & \multirow{3}{*}{$\mathrm{SO}(2)$} & uniform & $2^{84}$ & & & & \\
    $3$ &  & uniform & $2^{21}$ & 35\,s & 1.2\,s & 270\,ms & 37\,s \\
    \bfseries\boldmath $1$ &  & \bfseries\boldmath uniform & \bfseries\boldmath $2^7$ & \bfseries\boldmath 2.1\,ms & \bfseries\boldmath $<$1\,ms & \bfseries\boldmath $<$1\,ms & \bfseries\boldmath 2.5\,ms \\
  \end{tabular}
\end{table}

\FloatBarrier

\end{document}